\documentclass{isipta2025}

\renewbibmacro*{name:family}[4]{%
  \ifuseprefix
    {\usebibmacro{name:delim}{#3#1}%
     \usebibmacro{name:hook}{#3#1}%
     \mkbibcompletenamefamily{%
       \ifdefvoid{#3}
         {}
         {\mkbibnameprefix{\MakeCapital{#3}}\isdot%
          \ifprefchar{}{\bibnamedelimc}}%
       \mkbibnamefamily{#1}\isdot}}
    {\usebibmacro{name:delim}{#1}%
     \usebibmacro{name:hook}{#1}%
     \mkbibcompletenamefamily{%
       \mkbibnamefamily{#1}\isdot}}}%

\renewbibmacro*{name:family-given}[4]{%
  \ifuseprefix
    {\usebibmacro{name:delim}{#3#1}%
     \usebibmacro{name:hook}{#3#1}%
     \mkbibcompletenamefamilygiven{%
       \ifdefvoid{#3}
         {}
         {\mkbibnameprefix{\MakeCapital{#3}}\isdot
          \ifprefchar{}{\bibnamedelimc}}%
       \mkbibnamefamily{#1}\isdot
       \ifdefvoid{#4}
         {}
         {\bibnamedelimd\mkbibnamesuffix{#4}\isdot}%
       \ifdefvoid{#2}
         {}
         {\revsdnamepunct\bibnamedelimd\mkbibnamegiven{#2}\isdot}}}
    {\usebibmacro{name:delim}{#1}%
     \usebibmacro{name:hook}{#1}%
     \mkbibcompletenamefamilygiven{%
       \mkbibnamefamily{#1}\isdot
       \ifdefvoid{#4}
         {}
         {\bibnamedelimd\mkbibnamesuffix{#4}\isdot}%
       \ifboolexpe{%
         test {\ifdefvoid{#2}}
         and
         test {\ifdefvoid{#3}}}
         {}
         {\revsdnamepunct}%
       \ifdefvoid{#2}
         {}
         {\bibnamedelimd\mkbibnamegiven{#2}\isdot}%
       \ifdefvoid{#3}
         {}
         {\bibnamedelimd\mkbibnameprefix{#3}\isdot}}}}

\usepackage{tikz} 
\usetikzlibrary{automata,positioning,arrows.meta,shapes,calc,matrix,backgrounds}
\usepackage{booktabs} 
\usepackage{subfig}
\usepackage{array} 
\usepackage{tabularx} 
\usepackage{nicefrac}
\usepackage{aligned-overset}

\newif\ifarxiv
\arxivtrue 

\newcommand{\stsp}{\mathcal{X}}
\newcommand{\rfns}[1][\stsp]{\reals^{#1}}
\newcommand{\utranop}{\overline{T}{}}
\newcommand{\ltranop}{\underline{T}{}}
\newcommand{\restrutranop}[1][\class]{\utranop_{#1}}
\newcommand{\restrutranopsmash}[1][\class]{\utranop_{\smash{#1}}}
\newcommand{\restrltranop}[1][\class]{\ltranop_{#1}}

\newcommand{\nonabsrestrutranop}[1][\ell-1]{\restrutranop[\classtnl^{\smash{#1}}]}
\newcommand{\nonabsrestrltranop}[1][\ell-1]{\restrltranop[\classtnl^{\smash{#1}}]}
\newcommand{\tranop}{T}

\newcommand{\pmfset}[1]{\mathcal{P}^{#1}}
\newcommand{\credset}[1]{\mathcal{Q}^{#1}}
\newcommand{\restrpmfset}[2][\class]{\mathcal{P}_{#1}^{#2}}
\newcommand{\restrcredset}[2][\class]{\mathcal{Q}_{#1}^{#2}}

\newcommand{\uprev}{\overline{\mathrm{E}}}
\newcommand{\lprev}{\underline{\mathrm{E}}}
\newcommand{\prev}{\mathrm{E}}
\newcommand{\pmf}{p}
\newcommand{\altpmf}{q}

\DeclarePairedDelimiterXPP{\restr}[2]{}{.}{\downarrow}{_{#2}}{#1}
\DeclarePairedDelimiterXPP{\cylext}[2]{}{.}{\uparrow}{^{#2}}{#1}
\DeclarePairedDelimiter{\set}{\{}{\}}

\DeclarePairedDelimiter{\abs}{\vert}{\vert}
\DeclarePairedDelimiter{\card}{\vert}{\vert}
\DeclarePairedDelimiter{\gr}{(}{)}

\DeclarePairedDelimiterX{\ooi}[2]{]}{[}{#1, #2}
\DeclarePairedDelimiterX{\cci}[2]{[}{]}{#1, #2}
\DeclarePairedDelimiterX{\inti}[2]{\{}{\}}{#1,\ldots,#2}

\newcommand{\class}{C}
\newcommand{\altclass}{D}
\newcommand{\uagraph}[1][\utranop]{\overline{\mathcal{G}}\gr{#1}}
\newcommand{\cgraph}[1][\utranop]{\mathcal{G}^{\mathrm{c}}\gr{#1}}

\newcommand{\classesmax}{\mathcal{X}_{\mathrm{m}}}
\newcommand{\classmaxk}[1]{\mathcal{X}_{\mathrm{m}, \smash{#1}}}
\newcommand{\classtl}{\mathcal{X}_{\mathrm{t}\blacktriangledown}}
\newcommand{\classtnl}{\mathcal{X}_{\mathrm{t}\triangledown}}
\newcommand{\classtlsmash}[1]{\mathcal{X}^{#1}_{\smash{\mathrm{t}\blacktriangledown}}}
\newcommand{\classtnlsmash}[1]{\mathcal{X}^{#1}_{\smash{\mathrm{t}\triangledown}}}

\newcommand{\limitset}[1][f]{\Omega_{#1}}
\newcommand{\limit}{\omega}
\newcommand{\per}[1][f]{p_{#1}}

\DeclareMathOperator{\ext}{ext}

\newcommand{\nats}{\mathbb{N}}
\newcommand{\nnegints}{\mathbb{Z}_{\geq0}}
\newcommand{\reals}{\mathbb{R}}
\newcommand{\nnegreals}{\mathbb{R}_{\geq0}}
\newcommand{\posreals}{\mathbb{R}_{>0}}
\newcommand{\nposreals}{\mathbb{R}_{\leq0}}

\newcommand{\indica}[1]{\mathbb{1}_{#1}}

\newcommand{\sta}{\mathsf{a}}
\newcommand{\stb}{\mathsf{b}}
\newcommand{\stc}{\mathsf{c}}
\newcommand{\std}{\mathsf{d}}
\newcommand{\ste}{\mathsf{e}}

\newcommand{\hyp}[2]{\ensuremath{\mathrm{H}_{#1}^{#2}}}
\newcommand{\chyp}[1]{\ensuremath{\mathrm{H}_{#1}}}

\newtheorem{runex}{Running example}
\Crefname{runex}{Running example}{Running examples}
\AtBeginEnvironment{runex}{%
  \pushQED{\qed}%
}
\AtEndEnvironment{runex}{\popQED\endexample}

\title{A convenient characterisation of convergent upper transition operators}

\author[1]{Jasper De~Bock}
\author[1]{Alexander Erreygers}
\author[2]{Floris Persiau}

\affil[1]{Foundations Lab for imprecise probabilities\\Ghent University\\Belgium}
\affil[2]{Department of Philosophy\\Carnegie Mellon University\\United States of America}
\begin{document}

\maketitle


\begin{abstract}
	Motivated by its connection to the limit behaviour of imprecise Markov chains, we introduce and study the so-called convergence of upper transition operators: the condition that for any function, the orbit resulting from iterated application of this operator converges.
	In contrast, the existing notion of `ergodicity' requires convergence of the orbit to a constant.
	We derive a very general (and practically verifiable) sufficient condition for convergence in terms of accessibility and lower reachability, and prove that this sufficient condition is also necessary whenever (i) all transient states are absorbed or (ii) the upper transition operator is finitely generated.
\end{abstract}



\section{Introduction}

\emph{Imprecise Markov chains} \citep{2002Kozine,2009Skulj,2009deCooman,2016DeCooman,2018Denk,2014Hermans} model the uncertain temporal evolution of the state of finite-state discrete-time systems, and they do so in a more general manner than (`precise') Markov chains \citep{1960KemenySnell-Finite,1997Norris-Markov} by allowing for partial probability specifications.
Let \(\stsp\) be the finite set of possible states of the system under study.
A Markov chain models uncertain dynamics using transition probabilities: for every state~\(x\in\stsp\), one has to specify a transition probability mass function (pmf)~\(\pmf_x\), and then \(\pmf_x\gr{y}\) is the probability of transitioning from state~\(x\) in the current time step to state~\(y\) in the next one.
These give rise to a \emph{transition operator} (or \emph{kernel})~\(\tranop\colon\reals^{\stsp}\to\reals^{\stsp}\), which maps every function \(f\in\reals^\stsp\) to
\begin{equation*}
	\tranop f
	\colon \stsp\to\reals
	\colon x\mapsto \prev_{\pmf_x}\gr{f}.
\end{equation*}
This transition operator comes in handy when determining the expectation of a function~\(f\) of the state after \(n\)~time steps, given that the system started in state~\(x\):
\begin{equation*}
	\prev\gr[\big]{f\gr{X_n}\mid X_0=x}
	= \tranop^n f\gr{x}.
\end{equation*}

In contrast, an imprecise Markov chain models uncertain dynamics using a family \(\gr{\pmfset{x}}_{x\in\stsp}\) of candidate transition pmfs, with the interpretation that if up to time step~\(n\) the states were \(x_0\), \dots, \(x_n\), the uncertainty about the state~\(X_{n+1}\) in the next time step~\(n+1\) is accurately modelled by \emph{some} pmf in~\(\pmfset{x_n}\), which, depending on the adopted interpretation, may either depend on \(x_0\), \dots, \(x_n\) or only on \(x_n\).
Both interpretations give rise to the same range of values for the conditional expectation of a function~\(f\) of the state after \(n\)~time steps~\citep[Theorem~11.4]{2014Hermans}, and one is then typically interested in the supremum $\smash{\uprev\gr[\big]{f\gr{X_n}\mid X_0=x}}$ or infimum $\smash{\lprev\gr[\big]{f\gr{X_n}\mid X_0=x}}$ of this range, called the \emph{upper and lower expectation}, with
\begin{align*}
	\uprev\gr[\big]{f\gr{X_n}\mid X_0=x}
	&= \utranop^n f\gr{x},\\
&=-\lprev\gr[\big]{-f\gr{X_n}\mid X_0=x},
\end{align*}
where \(\utranop\) is the chain's \emph{upper transition operator}, which maps every function~\(f\in\rfns\) to
\begin{equation}
\label{eqn:utranop}
	\utranop f
	\colon \stsp\to\reals
	\colon x\mapsto \sup\set{\prev_{\pmf}\gr{f} \colon \pmf\in\pmfset{x}}.
\end{equation}
Studying the limit behaviour of imprecise Markov chains is therefore a matter of understanding the limit behaviour of the orbit~$\gr{\utranop^n f}_{n\in\nats}$.
In this contribution, we focus on the convergence of these orbits and develop conditions that guarantee convergence for all $f\in\rfns$.

\section{Upper transition operators}
Since in \cref{eqn:utranop} the upper transition operator~\(\utranop\) is defined as a pointwise supremum over linear operators, it has the following properties:
\begin{enumerate}[label=\upshape{T\arabic*.}, ref=\upshape T\arabic*, series=uto]
	\item\label{def:utranop:subadditive} \(\utranop \gr{f+g}\leq\utranop f+\utranop g\) for all \(f,g\in\rfns\);
	\item\label{def:utranop:positively homogeneous} \(\utranop \gr{\lambda f}=\lambda\utranop f\) for all \(f\in\rfns\) and \(\lambda\in\nnegreals\);
	\item\label{def:utranop:dominated} \(\utranop f\leq\max f\) for all \(f\in\rfns\).
\end{enumerate}
In general, we call any operator~\(\utranop\colon\rfns\to\rfns\) that satisfies \ref{def:utranop:subadditive}--\ref{def:utranop:dominated} an \emph{upper transition operator},\footnote{
	Also known as a \emph{sublinear transition operator} or a \emph{sublinear kernel} \citep[Definition~5.1]{2018Denk}.
} as it can always be thought of as being derived from sets of pmfs: \ref{def:utranop:subadditive}--\ref{def:utranop:dominated} guarantee that for all \(x\in\stsp\), there is at least one \emph{compatible} set of pmfs~\(\pmfset{x}\), in the sense that it satisfies the assignment in \cref{eqn:utranop}.
There is also a largest one, and this so-called \emph{credal set}~\(\credset{x}\) is the unique closed and convex set of pmfs that is compatible with $\utranop$. 



\subsection{Convergence and ergodicity} As explained in the Introduction, we aim to obtain (necessary and/or sufficient) conditions on the upper transition operator~\(\utranop\) for it to have the following property.
\begin{definition}
\label{def:convergent}
	An upper transition operator~\(\utranop\) is \emph{convergent} if for all \(f\in\rfns\), \(\smash{\gr{\utranop^n f}_{n\in\nats}}\) converges.
\end{definition}

Although we are not the first to study the limit behaviour of orbits~\(\smash{\gr{\utranop^n f}_{n\in\nats}}\), we are not aware of any work that focuses on (what we call) convergence.
Most work---at least on upper transition operators for imprecise Markov chains---focuses on the stronger requirement of ergodicity \cite{2009deCooman,2012Hermans,2011Skulj}; \citet{2013Skulj} does consider non-ergodic imprecise Markov chains, but focusses on invariant sets of distributions.
\begin{definition}
\label{def:ergodicity}
	An upper transition operator~\(\utranop\) is \emph{ergodic} if for all \(f\in\rfns\), \(\smash{\gr{\utranop^n f}_{n\in\nats}}\) converges to a constant.
\end{definition}
\Citet[Proposition~3]{2012Hermans} study this notion of ergodicity extensively, and they obtain a necessary and sufficient condition that is easy to check, which we'll introduce in \cref{sec:graph-theoretic} further on.
We aim for similar conditions, but for convergence rather than ergodicity.

\Citet{2012Hermans} also explain that thanks to \ref{def:utranop:subadditive}--\ref{def:utranop:dominated}, an upper transition operator~\(\smash{\utranop}\) is a `(convex) topical map' or, more generally, a `sup-norm non-expansive map' \citep{2005Lemmens,1990Sine}.
For these more general types of maps, the limit behaviour of orbits has been studied extensively as well \citep{1990Sine,2005Lemmens,2003Akian}.
From these references, we know that for any \(f\in\rfns\), the orbit~\(\gr{\utranop^n f}_{n\in\nats}\) has a finite limit set~\(\limitset=\set{\limit_1, \dots, \limit_{\per}}\),\footnote{
	Defined as the set of accumulation points of the orbit~\(\gr{\utranop^n f}_{n\in\nats}\), or equivalently, the limits of the convergent subsequences.
} whose \emph{period} (or cardinality)~\(\per\) has a universal upper bound that depends only on the size of~\(\stsp\) and whose elements form a cycle: \(\limit_2=\utranop\limit_1\), \dots, \(\limit_{\per}=\utranop\limit_{\per-1}\) and \(\limit_{\per+1}\coloneq\limit_1=\utranop\limit_{\per}\).
Convergence requires that $\per=1$ for all \(f\in\rfns\).

\subsection{Convenient properties of upper transition operators}
Throughout this contribution, we'll make use of several properties of upper transition operators that are well known \cite{2012Hermans,2012Hermans,2021TJoens}, but which we repeat here for the sake of convenience.
An upper transition operator~\(\utranop\) has a corresponding conjugate \emph{lower transition operator}~\(\ltranop\colon\rfns\to\rfns\) defined by~\(\ltranop f\coloneq-\utranop(-f)\) for all~\(f\in\rfns\); these names make sense because \(\ltranop f\leq \utranop f\) for all \(f\in\rfns\).
We'll often implicitly use that if \(\utranop\) is an upper transition operator, then so is its \(n\)-fold composition~\(\utranop^n\)---which has \(\ltranop^n\) as conjugate lower transition operator.
Other important properties of \(\smash{S\in\set{\utranop, \ltranop}}\) are:
\begin{enumerate}[resume*]
	\item\label{def:utranop:bounds} \(\min f\leq S f\leq\max f\) for all \(f\in\rfns\);
	\item\label{def:utranop:monotonicity} if \(f\leq g\) then \(S f\leq S g\) for all \(f,g\in\rfns\);
	\item\label{def:utranop:constant additivity} \(S \gr{\mu+f}=\mu+S f\) for all \(f\in\rfns\) and \(\mu\in\reals\);
	\item\label{prop:utranop:powerhorse} \((\max f-\min f)S \indica{x}+\min f \leq S f\) for all \(f\in\rfns\) and \(x\in\mathrm{arg\,max} f\).
\end{enumerate}
One particular type of functions we'll use are indicators: the \emph{indicator}~\(\indica{A}\) of the set~\(A\subseteq\stsp\) maps \(x\in\stsp\) to~\(1\) if \(x\in A\) and to~\(0\) otherwise; to ease our notation, we write \(\indica{x}\coloneq\indica{\set{x}}\) for all \(x\in\stsp\).
The final property we list revolves around indicators:
\begin{enumerate}[resume*]
	\item \label{prop:utranop:indicators} \(\utranop\indica{\class}=1-\ltranop\indica{\stsp\setminus\class}\) for all \(\class\subseteq\stsp\).
\end{enumerate}

\section{Characterising ergodicity}
\label{sec:graph-theoretic}

The necessary and sufficient condition for ergodicity given by \citet[Proposition~3]{2012Hermans} is stated in terms of accessibility and reachability relations.
Since these relations will also be important in our quest for conditions for convergence, we'll briefly recall them before repeating their result.
Throughout this contribution, we'll illustrate most relations and notions with the following running example.
\begin{runex}
	\label{example:running example first part}
	Let \(\stsp\coloneq\set{\sta, \stb, \stc,\std,\ste}\) and consider the upper transition operator~\(\utranop\) induced by the sets
\begin{multline*}
	\pmfset{\sta}\coloneq\set*{\indica{\sta}}
	\text{, }
	\pmfset{\stb}\coloneq\set*{\indica{\stb}}
	\text{, }
	\smash[t]{\pmfset{\stc}\coloneq\set*{\frac14\gr{\indica{\sta}+\indica{\stb}+\indica{\std}+\indica{\ste}}}} \\
	\text{and }
	\pmfset{\std}\coloneq\set*{\indica{\stc},\indica{\std},\indica{\ste}}\eqcolon\pmfset{\ste}.\qedhere
	\end{multline*}
\end{runex}

\subsection{The upper accessibility graph}
An upper transition operator~\(\utranop\) gives rise to a corresponding \emph{upper accessibility graph}~\(\uagraph\), with \(\stsp\) as nodes and a directed edge between states~\(x,y\in\stsp\) if and only if \(\utranop\indica{y}\gr{x}>0\) \citep[Definition~5]{2012Hermans}.
A state~\(y\in\stsp\) is now said to be \emph{accessible} from a state~\(x\in\stsp\) if either \(y=x\) or there is a directed path from \(x\) to~\(y\) in~\(\uagraph\), or equivalently, if \(\utranop^n\indica{y}\gr{x}>0\) for some \(n\in\nnegints\) \citep[Definition~4 and Proposition~4]{2012Hermans}.
\begin{runex}
	\label{example:running example first part bis}
	The reader will have no difficulty in verifying that \cref{fig:uagraph} depicts the upper accessibility graph~\(\uagraph\).
	\begin{figure}
	\centering
	\begin{tikzpicture}[node distance=1.5cm, on grid, auto]
		\node[state] (c) {\(\stc\)};
		\node[state] (a) [above left=of c] {\(\sta\)};
		\node[state] (b) [above right=of c] {\(\stb\)};
		\node[state] (d) [below left=of c] {\(\std\)};
		\node[state] (e) [below right=of c] {\(\ste\)};
		\path[->,>={Stealth[round]}]
		(a) edge[loop left]  (a)
		(b) edge[loop right] (b)
		(d) edge[loop left]  (d)
		(e) edge[loop right] (e)
		(c) edge             (a)
		(c) edge             (b)
		(c) edge[bend right] (d)
		(c) edge[bend left]  (e)
		(d) edge             (c)
		(d) edge[bend right] (e)
		(e) edge             (c)
		(e) edge[bend right] (d);
	\end{tikzpicture}
	\caption{\(\uagraph\) for \cref{example:running example first part}.}
	\label{fig:uagraph}
	\end{figure}
\end{runex}

Two states~\(x,y\in\stsp\) \emph{communicate} if \(y\) is accessible from~\(x\) (in~\(\uagraph\)) and vice versa; this equivalence relation partitions the state space~\(\stsp\) into equivalence classes \(\class_1\), \dots, \(\class_n\), aptly called \emph{communication classes}---or strongly connected components in the theory of directed graphs.
For two communication classes $C_k$ and $C_\ell$, $y\in C_\ell$ is accessible from $x\in C_k$ if and only if the same is true for any $y'\in C_\ell$ and $x'\in C_k$; whenever this is the case, we therefore simply say that $C_\ell$ is accessible from~$C_k$.
This accessibility relation induces a partial order on the communication classes: we say that $C_\ell$ dominates $C_k$ if $C_\ell\neq C_k$ and $C_\ell$ is accessible from $C_k$.
A communication class~\(\class\) is called \emph{maximal} (or final) if it's undominated with respect to this partial order. Since there are only finitely many communication classes, every non-maximal communication class is dominated by at least one maximal communication class. A state is called \emph{maximal} if it belongs to a maximal communication class and \emph{non-maximal} or \emph{transient} otherwise.
We enumerate the maximal communication classes by \(\classmaxk{1}\), \dots, \(\classmaxk{M}\), and we collect their union in~\(\classesmax\coloneq\classmaxk{1}\cup\cdots\cup\classmaxk{M}\).

The \emph{cyclicity} or \emph{period} of a maximal communication class~\(\classmaxk{k}\) is the greatest common divisor of the lengths of the closed directed paths that remain in this class, and the class~\(\classmaxk{k}\) is said to be \emph{regular} if it has cyclicity~\(1\), or equivalently, if there is some \(N\in\nats\) such that for all~\(n\geq N\) and \(x,y\in\classmaxk{k}\), \(\utranop^n\indica{y}\gr{x}>0\) \citep[Proposition~4.2]{2009deCooman}.

\begin{runex}
\label{example:running example second part}
	The graph~\(\smash{\uagraph}\) has three communication classes: \(\set*{\sta}\), \(\set*{\stb}\) and \(\set*{\stc,\std,\ste}\).
	The first two are the maximal ones, and these are regular because they obviously have cyclicity~\(1\).
\end{runex}

For small state spaces, it's easy to determine the maximal communication classes and their cyclicity on sight.
In general, however, there's several of algorithms that do this (in linear time)---see for example \citep[Sections~13.2.3 and 13.3.2]{1999Jarvis}.

To make our lives a bit easier, we'll call any non-empty subset~\(\class\) of~\(\stsp\) a \emph{class}.
Such a class~\(\class\) is \emph{closed}---after \citet[Section~5]{2021TJoens}---if \(\utranop\indica{y}\gr{x}=0\) for all \(x\in \class\) and \(y\in\stsp\setminus\class\); it follows from \ref{def:utranop:subadditive} and \ref{def:utranop:bounds} that this is the case if and only if \(\utranop\indica{\stsp\setminus\class}\gr{x}=0\) for all \(x\in\stsp\).
Clearly, a communication class~\(\class\) is maximal if and only if it's closed; we leave it to the reader to verify that any closed class must be a union of communication classes, but that not all such unions need be closed.

\begin{runex}
\label{example:running example second part extra}
	The closed classes in our running example are \(\set*{\sta}\), \(\set*{\stb}\),  \(\set*{\sta,\stb}\) and \(\stsp\).
\end{runex}

\subsection{Lower reachability and absorption}
Unfortunately, the upper accessibility graph~\(\smash{\uagraph}\) does not suffice to characterise ergodicity; for this, we also need the notions of `lower reachability' and `absorption'.
A class~\(\class\) is \emph{lower reachable} from a state~\(x\in\stsp\) if there is some \(n\in\nats\) such that \(\ltranop^n\indica{\class}(x)>0\), and we call \(\class\) \emph{absorbing} if it is lower reachable from any state~\(x\in\stsp\setminus\class\) \citep[Section~5]{2021TJoens}.
The following result---a slight generalisation of \citep[Proposition~6]{2012Hermans}---provides a convenient recursive method to determine the states from which a \emph{closed} class is lower reachable.
\begin{lemma}\label{lem:iterativeabs}
	Consider an upper transition operator~\(\utranop\) with closed class \(\class\).
	Let \(\gr{C_n}_{n\in\nnegints}\) be the non-decreasing sequence given by \(\class_0\coloneq\class\) and, for all \(n\in\nnegints\), by
	\begin{align*}
	\class_{n+1}
	\coloneq{}&\class_n\cup\{x\in\stsp\setminus\class_n\colon \ltranop\indica{\class_n}(x)>0\}\\
	={}& \set{x\in\stsp\colon \ltranop^{n+1} \indica{\class}\gr{x}>0}.
	\end{align*}
	Then after \(k\leq\card{\stsp\setminus\class}\) steps, we reach \(\class_k=\class_{k+1}\), and \(\class_k\) is the set of states from which~\(\class\) is lower reachable.
\end{lemma}
\ifarxiv
\begin{proof}
Observe that \(\gr{\class_n}_{n\in\nats}\) is a non-decreasing sequence by definition, and that \(\stsp\setminus\class_0=\stsp\setminus\class\) is finite.
This implies that there must be some first \(k\in\nnegints\) such that \(\class_{k+1}=\class_k\)---and then \(\class_k=\class_n\) for all \(n\geq k\)---and that \(k\leq\card{\stsp\setminus\class_0}=\card{\stsp\setminus\class}\).
So to prove the statement, it remains for us to verify that
\begin{equation*}
	\class_n
	= \set{x\in\stsp\colon \ltranop^n \indica{\class}\gr{x}>0}
	\quad\text{for all } n\in\nnegints.
\end{equation*}
To this end, it suffices to verify, for all \(n\in\nnegints\), the veracity of the conjunction~\(\chyp{n}\coloneq\hyp{n}{1}\wedge\hyp{n}{2}\wedge\hyp{n}{3}\) of the following three statements:
\begin{align*}
	\hyp{n}{1}
	&\colon \gr{\forall x\in\stsp\setminus\class_n}~\ltranop^n\indica{\class}(x) = 0;
	\\
	\hyp{n}{2}
	&\colon\gr{\forall x\in \class_n}~\ltranop^n\indica{\class}(x)>0; \\
	\hyp{n}{3}
	&\colon\gr{\forall x\in \class_{n+1}}~\ltranop\indica{\class_n}(x)>0.
\end{align*}
Our proof that \(\chyp{n}\) holds for all \(n\in\nnegints\) is one by induction.
For the base case \(n=0\), \hyp{0}{1} and \hyp{0}{2} are trivially true because \(\class=\class_0\) by definition, while \hyp{0}{3} holds because (i) for all \(x\in\class_1\setminus\class_0\), \(\ltranop\indica{\class_0}\gr{x}>0\) by definition, and (ii) for all \(x\in\class_0\),
\begin{equation*}
	\ltranop \indica{\class_0}\gr{x}
	\overset{\smash[t]{\text{\ref{prop:utranop:indicators}}}}{=} 1-\utranop\indica{\stsp\setminus\class_0}\gr{x}
	= 1,
\end{equation*}
where for the last equality we used the assumption in the statement that \(\class_0=\class\) is a closed class.

For the inductive step, we fix some \(n\in\nnegints\) and show that \(\chyp{n}\Rightarrow\chyp{n+1}\).
It follows from \hyp{n}{1} and \hyp{n}{2} that, with \(\delta\coloneq\min_{x\in\class_n}\ltranop^n\indica{\class}(x)>0\), \(\indica{\class_n}\geq\ltranop^n\indica{\class}\geq\delta\indica{\class_n}\) and therefore
\begin{equation*}
	\ltranop \indica{\class_n}
	\overset{\smash[t]{\text{\ref{def:utranop:monotonicity}}}}{\geq} \ltranop^{n+1} \indica{\class}
	\overset{\smash[t]{\text{\ref{def:utranop:positively homogeneous},\ref{def:utranop:monotonicity}}}}{\geq} \delta\ltranop\indica{\class_n}.
\end{equation*}
From these inequalities, we infer that
\begin{align*}
	\set{x\in\stsp\colon \ltranop^{n+1}\indica{\class}\gr{x}>0}
	= \set{x\in\stsp\colon \ltranop\indica{\class_n}\gr{x}>0}.
\end{align*}
Recalling that \(\class_{n+1}=\class_n\cup\set{x\in\stsp\setminus\class_n\colon\ltranop\indica{\class_n}\gr{x}>0}\) by definition, it follows from \hyp{n}{3} that
\begin{equation*}
	\class_{n+1}
	= \set{x\in\stsp\colon \ltranop\indica{\class_n}\gr{x}>0}.
\end{equation*}
The preceding two equalities give us \hyp{n+1}{1} and \hyp{n+1}{2}.
Statement \hyp{n+1}{3} holds because (i) for all \(x\in\class_{n+2}\setminus\class_{n+1}\), \(\ltranop\indica{\class_{n+1}}\gr{x}>0\) by definition; and (ii) for all \(x\in\class_{n+1}\),
\begin{equation*}
	\ltranop\indica{\class_{n+1}}\gr{x}
	\overset{\smash[t]{\text{\ref{def:utranop:monotonicity}}}}{\geq} \ltranop\indica{\class_n}\gr{x}
	> 0,
\end{equation*}
where for the first inequality we also used that \(\class_{n+1}\supseteq\class_n\) by definition, and for the second inequality we used \hyp{n}{3}.
\end{proof}
\fi
As \(\classesmax\) is the union of the maximal (and therefore closed) communication classes, it is closed as well and we can determine the set $C_k$ of states from which it is lower reachable with the iterative method in~\cref{lem:iterativeabs} for \(\class=\classesmax\).
We collect all transient states of this kind in~\(\classtl=\class_k\setminus\classesmax\); finally, \(\classtnl\coloneq\stsp\setminus\gr{\classesmax\cup\classtl}\) collects the transient states from which~\(\classesmax\) is \emph{not} lower reachable.
Note that, as depicted in~\cref{fig:one cut}, \(\classtl\) and \(\classtnl\) need \emph{not} be communication classes, and that \(\gr{\classesmax, \classtl, \classtnl}\) partitions~\(\stsp\).
\begin{figure}
	\centering
	\begin{tikzpicture}[divh/.style={minimum height=2.4em, inner sep=2pt}, div/.style={minimum width=3.8em, divh}, partit/.style={draw, rounded corners, thick}, max/.style={div, partit, color=teal}, abson/.style={divh, color=magenta}, absob/.style={partit, color=magenta}, nonabson/.style={divh, color=cyan}, nonabsob/.style={partit, color=cyan}, commcl/.style={draw, rectangle, rounded corners, dashed},]
		\node[max, dashed] (max1) {\(\classmaxk{1}\)};
		\node[div, color=teal, right = .25cm of max1] (cdots) {\(\cdots\)};
		\node[max, dashed, right = .25cm of cdots] (cmaxM) {\(\classmaxk{M}\)};
		\node[divh, right = .25cm of cmaxM, color=teal] (cmax) {\(\classesmax\)};
		\draw[max] ($(max1.north west)+(-.25, .25)$) rectangle ($(cmax.south east)+(.25, -.25)$);

		\node[abson, below= .5cm of max1, xshift=-.25cm] (ca) {\(\classtl\)};
		\node[nonabson, below = .5cm of cmax] (cna) {\(\classtnl\)};

		\coordinate[xshift=-.25cm] (helperleft) at (max1.south west|-ca.south west);
		\coordinate (helpermid) at (cna.north west-|cmaxM.south west);
		\coordinate[xshift=.25cm] (helperright) at (cna.south east-|cmax.south east);
		\draw[absob] (helperleft) rectangle ($(helpermid)+(-.25, 0)$);
		\draw[nonabsob]	(helpermid) rectangle (helperright);

		\draw[commcl] ($(ca.north east)+(.25, -.125)$) rectangle ($(ca.south east)+(1, .125)$);
		\draw[commcl] ($(cna.north west)+(-.25, -.125)$) rectangle ($(cna.south west)+(-.75, .125)$);
		\draw[commcl] ($(ca.north east)+(1.25, -.125)$) rectangle ($(cna.south west)+(-1, .125)$);

		\draw[partit] ($(max1.north west)+(-.5, .5)$) rectangle ($(helperright)+(.25, -.25)$);
	\end{tikzpicture}
	\caption{Venn diagram of the state space~\(\stsp\)}
	\label{fig:one cut}
\end{figure}
\begin{runex}
\label{example:running example third part}
	To determine which transient states \(\classesmax\) is lower reachable from, we apply the recursive method in \cref{lem:iterativeabs} for \(\class=\classesmax=\set{\sta, \stb}\):
	\begin{align*}
	\class_1
	&=
	\set{\sta,\stb}\cup\set{x\in\set{\stc,\std,\ste}\colon\ltranop\indica{\set{\sta,\stb}}(x)>0}
	=
	\set{\sta,\stb,\stc}
	\shortintertext{and}
	\class_2
	&=
	\set{\sta,\stb,\stc}\cup\set{x\in\set{\std,\ste}\colon\ltranop\indica{\set{\sta,\stb,\stc}}(x)>0}
	=
	\class_1.
	\end{align*}
	Consequently, \(\classesmax\) is lower reachable from \(\stc\) but not from \(\std\) and \(\ste\).
	This illustrates nicely that \(\classtl=\set{\stc}\) and \(\classtnl=\set{\std,\ste}\) need not be communication classes.
\end{runex}

Crucially, the class \(\classtnl\) of transient states from which \(\classesmax\) is not lower reachable is always empty for a (linear) transition operator~\(\tranop\).
\begin{lemma}
\label{lem:linear then classesmax absorbing}
    For any linear transition operator~\(\tranop\), \(\classesmax\) is absorbing, so \(\classtl=\stsp\setminus\classesmax\) and \(\classtnl=\emptyset\).
\end{lemma}
\begin{proof}
    The statement is trivial if \(\classesmax=\stsp\), so assume that \(\stsp\setminus\classesmax\neq\emptyset\).
    Any non-maximal \(x\in\stsp\setminus\classesmax\) belongs to a non-maximal communication class \(\class\), for which we know that it is dominated by some maximal communication class \(\classmaxk{k}\).
    For any $y\in\classmaxk{k}\subseteq\classesmax$, this implies that $y$ is accessible from $x$, meaning that there is some \(n\in\nats\) such that \(\tranop^n \indica{y}\gr{x}>0\) and therefore also \(\tranop^n \indica{\classesmax}\gr{x}\geq\tranop^n \indica{y}\gr{x}>0\).
    Since $\tranop$ is linear, this shows that \(\classesmax\) is lower reachable from any \(x\in\stsp\setminus\classesmax\), whence \(\classesmax\) is absorbing.
\end{proof}

\subsection{A necessary and sufficient condition for ergodicity}
And with that, we are ready to state the necessary and sufficient condition for ergodicity given by \citet[Proposition~3]{2012Hermans}:
\begin{proposition}
\label{prop:ergodic iff}
	An upper transition operator~\(\utranop\) is ergodic if and only if it has a single maximal communication class (so \(\classesmax=\classmaxk{1}\)) that is absorbing (so \(\classtnl=\emptyset\)) and regular.
\end{proposition}
Since ergodicity clearly implies convergence, this result provides a sufficient condition for convergence that can be checked easily.
However, whenever there's more than one maximal communication class, as in our running example, we're already out of luck.
\begin{runex}
	\label{example:running example fourth part}
	Since the upper transition operator~\(\utranop\) has two maximal communication classes, \cref{prop:ergodic iff} tells us it cannot be ergodic.
	However, it \emph{is} convergent!
	The reader may set out to verify this explicitly, but it is much more convenient to use the sufficient condition for convergence that we'll establish in \cref{the: sufficient condition for convergence} further on---see \cref{example:running example fifth part}.
\end{runex}

\section{Conditions for convergence}
\label{sec:conditions for convergence}
Our path forward is clear: we set out to establish a sufficient condition for convergence that is more general than the (necessary and sufficient) one for ergodicity in \cref{prop:ergodic iff}.
We'll do so in several stages: we first obtain a necessary and sufficient condition in the case of a single communication class, then a necessary and sufficient condition when there's multiple maximal communication classes whose union is absorbing, and finally a sufficient condition in the general case.

\subsection{A single communication class}
Let us start gently with upper transition operators that have a single communication class---a linear transition operator of this type is called \emph{irreducible} \citep[Section~13.2.1]{1999Jarvis}.
We know from \cref{prop:ergodic iff} that in this case regularity is sufficient for ergodicity and therefore convergence, and the following result establishes that it is also necessary for convergence.
\begin{proposition}\label{prop:equivalence single communication class}
	Consider an upper transition operator~\(\utranop\) with a single communication class.
	Then the following three statements are equivalent: \emph{(i)} $\utranop$ is convergent; \emph{(ii)} $\utranop$ is ergodic; \emph{(iii)} $\stsp$ is regular.
\end{proposition}
Besides \cref{prop:ergodic iff}, our proof for this result relies on the following lemma, which will come in handy further on as well.
\begin{lemma}\label{lem:limitiffconstant}
Consider an upper transition operator~\(\utranop\) with a single communication class.
Then for any \(f\in\rfns\), \(\gr{\utranop^n f}_{n\in\nats}\) converges if and only if it converges to a constant.
\end{lemma}
\begin{proof}
It suffices to show that if \(\phi\coloneq\lim_{n\to+\infty}\utranop^n f\) exists, it must be constant.
So let us assume \emph{ex absurdo} that $\phi$ exists but is not constant.
Then there are $x,y\in\stsp$ such that $\phi(x)=\max\phi>\min\phi=\phi(y)$.
Since $\stsp$ is a communication class, there is some $k\in\nats$ such that $\utranop^k\indica{x}(y)>0$.
It therefore follows that
\begin{equation*}
\utranop^k\phi(y)
\overset{\smash[t]{\text{\ref{prop:utranop:powerhorse}}}}{\geq}(\max \phi-\min \phi)\utranop^k\indica{x}(y)+\min \phi>\min\phi.
\end{equation*}
Meanwhile, \(\utranop\phi=\phi\), and hence, \(\phi(y)=\utranop^k\phi(y)>\min\phi\), which contradicts the fact that \(\phi(y)=\min\phi\).
\end{proof}
Proving \cref{prop:equivalence single communication class} is now a piece of cake.
\begin{proof}[Proof of \cref{prop:equivalence single communication class}]
The equivalence of convergence and ergodicity is an immediate consequence of \cref{lem:limitiffconstant}. That ergodicity is equivalent to regularity, on the other hand, follows directly from \cref{prop:ergodic iff}
\end{proof}

Whenever
one of the three equivalent conditions in \cref{prop:equivalence single communication class} holds, we can also say something about the (constant) limit of~\(\gr{\utranop^n f}_{n\in\nats}\) in relation to~\(f\).
\begin{proposition}\label{prop:limitdominates}
Consider an upper transition operator~\(\utranop\) such that \(\stsp\) is a regular communication class.
Then for any \(f\in\rfns\), \(\gr{\utranop^n f}_{n\in\nats}\) converges to a constant function~$\phi\geq\min f$ and, unless $f$ is constant, $\phi>\min f$.
\end{proposition}
\begin{proof}
Since $\stsp$ is a regular communication class, we know from \cref{prop:equivalence single communication class} that $\utranop$ is ergodic, which implies that \(\gr{\utranop^n f}_{n\in\nats}\) converges to a constant function $\phi$.

Consider any $x\in\mathrm{arg\,max} f$. Since $\stsp$ is regular, there is some $k\in\nats$ such that $\min\utranop^k\indica{x}>0$.
Consequently, with \(\alpha\coloneq\max f-\min f\geq0\),
\begin{equation*}
\phi=\lim_{n\to+\infty}\utranop^n\utranop^kf
\overset{\smash[t]{\text{\ref{def:utranop:bounds}}}}{\geq}\min\utranop^kf
\overset{\smash[t]{\text{\ref{prop:utranop:powerhorse}}}}{\geq}\alpha\min\utranop^k\indica{x}+\min f.
\end{equation*}
Hence, indeed, $\phi\geq\min f$; if $f$ is not constant, then \(\alpha>0\) and therefore $\phi>\min f$.
\end{proof}

\subsection{Restriction}
In the remainder of this section we move beyond upper transition operators with a single communication class. As a first step, instead of studying the convergence of \(\utranop\) for all $x\in\stsp$, we zoom in on particular subsets of states.

To do so, we introduce the notion of restriction of functions, sets of pmfs and upper transition operators.
For \(A,B\subseteq\stsp\) with \(A\supseteq B\) and \(f\in\rfns[A]\), \(\restr{f}{B}\) denotes the \emph{restriction of~\(f\) to~\(B\)}:
\begin{equation*}
	\restr{f}{B}
	\colon B\to\reals
	\colon x\mapsto f(x).
\end{equation*}
For any family~\(\gr{\pmfset{x}}_{x\in\stsp}\) of sets of pmfs, class \(\class\subseteq\stsp\) and \(x\in\stsp\), we let
\begin{equation}
	\restrpmfset{x}
	\coloneq \set{\restr{\pmf}{\class}\colon \pmf\in\pmfset{x}, (\forall y\in\stsp\setminus \class)~\pmf\gr{y}=0}
	\label{eqn:restrcredset:masses}
\end{equation}
be the set that collects the restriction to the class~\(\class\) of those pmfs~\(\pmf\in\pmfset{x}\) whose support~\(\set{y\in\stsp\colon \pmf\gr{y}>0}\) is contained in~\(\class\).
Consider now an upper transition operator~\(\smash{\utranop}\) with corresponding family~\((\credset{x})_{x\in\stsp}\) of credal sets. 
For any class~\(\class\subseteq\stsp\) such that \(\restrcredset{x}\neq\emptyset\) for all \(x\in \class\), we let \(\smash{\restrutranop}\) be the upper transition operator corresponding to the family~\((\restrcredset{x})_{x\in\class}\):
\begin{equation*}
	\restrutranop g\gr{x}
	= \max\set{\prev_\pmf(g)\colon \pmf\in\restrcredset{x}}
	\;\text{for all } g\in\rfns[\class], x\in\class.
\end{equation*}
It follows almost immediately from the definition of~\(\restrutranop\) and \(\gr{\restrcredset{x}}_{x\in\stsp}\) that
\begin{equation}
\label{eqn:inequality with restrutranop}
	\gr{\restrutranopsmash}^n (\restr{f}{\class})
	\leq \restr{(\utranop^n f)}{\class}
	\quad\text{for all } n\in\nats, f\in\rfns.
\end{equation}
\begin{proof}
	The base case \(n=1\) follows immediately from the definition of~\(\restrutranop\) and \(\restrpmfset{x}\).
	The inductive step follows from this base case and \ref{def:utranop:monotonicity}.
\end{proof}

One thing that is particularly useful in practice, is that one can obtain \(\restrutranop\) with \emph{any} compatible family of sets of pmfs, as long as these sets are closed.
\begin{lemma}
\label{lem:restrutranop via closed family}
	Consider an upper transition operator~\(\utranop\), let \(\gr{\pmfset{x}}_{x\in\stsp}\) be any compatible family of \emph{closed} sets of pmfs and fix some class \(\class\subseteq\stsp\).
	Then \(\restrutranop\) is well defined if and only if for all \(x\in\class\), \(\utranop\indica{\class}\gr{x}=1\), or equivalently, \(\restrpmfset{x}\neq\emptyset\); whenever this is the case, \(\smash{\gr{\restrpmfset{x}}_{x\in\class}}\) is compatible with~\(\smash{\restrutranop}\).
\end{lemma}
\ifarxiv
Our proof for this result requires a little bit extra information.
For any closed and convex set~\(S\), we denote the set of extreme points---those pmfs in \(S\) that cannot be written as convex combinations of others \citep[Section~18]{1970Rockafellar-Convex}---by~\(\ext S\).
Furthermore, we'll rely on the well-known result that for any upper transition operator~\(\utranop\), \(\gr{\ext \credset{x}}_{x\in\stsp}\) is compatible with~\(\utranop\)---see for example \citep[Proposition~2.5]{2014DeCoomanMiranda}.
\begin{proof}[Proof of \cref{lem:restrutranop via closed family}]
    For the first part of the statement, we show that for all \(x\in\class\),
    \begin{equation*}
        \utranop \indica{\class}\gr{x}
        = 1
        \Rightarrow \restrpmfset{x}\neq\emptyset
        \Rightarrow \restrcredset{x}\neq\emptyset
        \Rightarrow \utranop \indica{\class}\gr{x}
        = 1.
    \end{equation*}
    The first implication follows because \(\pmfset{x}\) is assumed to be closed: then there is some \(\pmf\in\pmfset{x}\) such that \(\prev_\pmf\gr{\indica{\class}}=1\), whence it must be that \(\pmf\gr{y}=0\) for all \(y\in\stsp\setminus\class\), so \(\restr{\pmf}{\class}\in\restrpmfset{x}\neq\emptyset\).
For the second implication, note that \(\pmfset{x}\subseteq\credset{x}\) because \(\gr{\pmfset{y}}_{y\in\stsp}\) is compatible with~\(\smash{\utranop}\) by assumption.
  From this and \cref{eqn:restrcredset:masses}, it follows immediately that \(\restrpmfset{x}\subseteq\restrcredset{x}\).
For the final implication, the left-hand side implies that there is some \(\pmf\in\credset{x}\) such that \(\pmf\gr{y}=0\) for all \(y\in\stsp\setminus\class\), whence \(\prev_\pmf\gr{\indica{\class}}=1\) and therefore \(\utranop\indica{\class}\gr{x}=1\).

    For the remainder of the statement, observe that for any \(x\in\class\), \(\credset{x}\) is closed and convex so \(\restrcredset{x}\) is clearly also closed and convex.
    We set out to prove that
	\begin{equation*}
		\ext\restrcredset{x}
		\subseteq \restrpmfset{x}
		\subseteq \restrcredset{x}
		\quad\text{for all } x\in\class.
	\end{equation*}
  Since \(\smash{\gr{\ext\restrcredset{x}}_{x\in\class}}\) and  \(\smash{\gr{\restrcredset{x}}_{x\in\class}}\) are both compatible with~\(\restrutranop\), this implies the second part of the statement.

	Fix any \(x\in\class\).
	We've already argued in the first part of this proof that \(\restrpmfset{x}\subseteq\restrcredset{x}\).

	The remaining inclusion \(\ext\restrcredset{x}\subseteq\restrpmfset{x}\) is clearly valid if \(\ext\restrcredset{x}=\emptyset\), so we assume that \(\ext\restrcredset{x}\neq\emptyset\) and fix some $q\in\ext\restrcredset{x}$. Then there is some \(\pmf\in\credset{x}\) such that \(\restr{\pmf}{\class}=q\).
	We now show that \(\pmf\) is an extreme point of~\(\credset{x}\).
	To that end, we assume towards contradiction that it is not.
	Then \(\pmf=\lambda\pmf_1+\gr{1-\lambda}\pmf_2\) for pmfs~\(\pmf_1\neq\pmf_2\) in~\(\credset{x}\) and \(\lambda\in\ooi{0}{1}\), but then \(\restr{\pmf_1}{\class}, \restr{\pmf_2}{\class}\in\restrcredset{x}\) with \(\restr{\pmf_1}{\class}\neq\restr{\pmf_2}{\class}\) and \(\restr{\pmf}{\class}=\lambda\restr{\pmf_1}{\class}+\gr{1-\lambda}\restr{\pmf_2}{\class}\), a contradiction because \(\restr{\pmf}{\class}=q\) is an extreme point of~\(\restrcredset{x}\).
	Therefore, \(\pmf\) is indeed an extreme point of~\(\credset{x}\).

	Let $\mathrm{conv}(\pmfset{x})$ be the convex hull of $\pmfset{x}$.
	On the one hand, since $\pmfset{x}\subseteq\mathrm{conv}(\pmfset{x})\subseteq\credset{x}$ and $\pmfset{x}$ and $\credset{x}$ are both compatible with $\utranop$, we know that $\mathrm{conv}(\pmfset{x})$ is compatible with $\utranop$.
	On the other hand, since $\pmfset{x}$ is closed, $\mathrm{conv}(\pmfset{x})$ is closed and convex.
	This implies that $\credset{x}=\mathrm{conv}(\pmfset{x})$ and therefore that $\ext\credset{x}=\ext\mathrm{conv}(\pmfset{x})$.
	Since \citep[Corollary 18.3.1]{1970Rockafellar-Convex} tells us that $\ext\mathrm{conv}(\pmfset{x})\subseteq\pmfset{x}$, we find that $\ext\credset{x}\subseteq\pmfset{x}$.
	Since \(\pmf\in\credset{x}\) is an extreme point of~\(\credset{x}\), we conclude that $\pmf\in\pmfset{x}$ and, consequently, that \(q=\restr{\pmf}{\class}\in\restrpmfset{x}\).
\end{proof}
\fi
\begin{runex}\label{example:running example eerste keer knippen}
Recall that \(\classesmax=\classmaxk{1}\cup\classmaxk{2}=\set{\sta,\stb}\), \(\classtl=\set{\stc}\) and \(\classtnl=\set{\std,\ste}\).
Applying \cref{eqn:restrcredset:masses} for $x\in\class=\classtnl$ we find that
\begin{equation*}
\restrpmfset[\classtnl]{\std}=\restrpmfset[\classtnl]{\ste}=\set{\indica{\std},\indica{\ste}}.
\end{equation*}
\cref{lem:restrutranop via closed family} therefore implies that \(\smash{\restrutranop[\classtnl]}\) is well defined---which is no coincidence, as we'll see in \cref{lem:nonabs has nonempty restrcredset}---and, in particular, that \(\restrutranop[\classtnl]f=\max f\) for all $\smash{f\in\rfns[\classtnl]}$; its upper accessibility graph is depicted in \cref{fig:running example second part}.
\begin{figure}
\centering
\begin{tikzpicture}[node distance=1.5cm, on grid, auto]
	\node[state] (d) {\(\std\)};
	\node[state] (e) [right=of d] {\(\ste\)};
	\path[->,>={Stealth[round]}]
	(d) edge[loop left]  (d)
	(e) edge[loop right] (e)
	(d) edge[bend right] (e)
	(e) edge[bend right] (d);
\end{tikzpicture}
\caption{\(\uagraph[{\restrutranop[\classtnl]}]\) for \cref{example:running example eerste keer knippen}.}
\label{fig:running example second part}
\end{figure}
\end{runex}

With these notions of restriction in place, we can now formalize what we mean by zooming in on $\class$.

\begin{definition}
\label{def:generalconvergence}
	An upper transition operator~\(\utranop\) is \emph{convergent} (\emph{ergodic}) \emph{on} $\class\subseteq\stsp$ if, for all \(f\in\rfns\), \(\gr{\restr{(\utranop^n f)}{\class}}_{n\in\nats}\) converges (to a constant).
\end{definition}

Whenever $\class$ is a maximal communication class, these restricted notions of convergence and/or ergodicity can be conveniently characterised in terms of~\(\restrutranop\).
\begin{lemma}\label{lem:closed classes}
Let \(\class\) be one of the maximal communication classes of an upper transition operator~\(\utranop\).
Then \(\restrutranop\) is well defined and
\begin{equation*}
		(\restrutranop)^n (\restr{f}{\class})
		= \restr{(\utranop^n f)}{\class}
		\quad\text{for all } n\in\nats, f\in\rfns.
\end{equation*}
Consequently,
\begin{enumerate}[label=\upshape(\roman*)]
\item \(\class\) is the sole communication class for~\(\restrutranop\);\label{item:one}
\item \(\class\) is regular for~\(\restrutranop\) if and only if it's regular for~\(\utranop\);\label{item:two}
\item $\restrutranop$ is ergodic if and only if $\utranop$ is ergodic on~$C$;\label{item:four}
\item $\restrutranop$ is convergent if and only if $\utranop$ is convergent on $C$.\label{item:five}
\end{enumerate}
\end{lemma}
\ifarxiv
\begin{proof}
	Let \(\gr{\credset{x}}_{x\in\stsp}\) be the family of credal sets corresponding to~\(\utranop\).
	Fix any \(x\in\class\).
	Then \(\utranop\indica{y}(x)=0\) for all \(y\in \stsp\setminus\class\) since the maximal communication class~\(\class\) is closed, and therefore
	\begin{equation}
	\label{eqn:pmfy0}
		\pmf(y)=0
		\quad\text{for all } \pmf\in\credset{x}, y\in\stsp\setminus\class.
	\end{equation}
	It follows from this and \cref{eqn:restrcredset:masses} that \(\credset{x}\neq\emptyset\).
	Furthermore, we infer from \cref{eqn:pmfy0} that for all \(f\in\rfns\),
	\begin{equation*}
		\utranop f(x)
		= \max_{\pmf\in\credset{x}} \prev_\pmf(f)
		= \max_{\altpmf\in\restrcredset{x}} \prev_\altpmf(\restr{f}{\class})
		= \restrutranop (\restr{f}{\class})(x).
	\end{equation*}
	This proves the equality in the statement for \(n=1\).
	The general case follows by induction: for all \(n\geq 2\), \(f\in\rfns\) and \(x\in\class\),
	\begin{align*}
		(\restrutranop)^n \gr{\restr{f}{\class}}(x)
		&= \restrutranop \gr[\big]{\gr{\restrutranop}^{n-1}\gr{\restr{f}{\class}}} (x) \\
		&= \restrutranop \gr[\big]{\restr{\gr{\utranop^{n-1}f}}{\class}} (x)
		= \utranop^n f(x).
	\end{align*}

	For statements~\ref{item:one} and \ref{item:two}, it suffices to realise that due to the first part of the statement,
	\begin{equation*}
		\utranop^n \indica{y}\gr{x}
		= (\restrutranop)^n \indica{y}\gr{x}
		\quad\text{for all } n\in\nats, x,y \in\class.
	\end{equation*}
	Statements~\ref{item:four} and \ref{item:five} follow immediately from the first part of the statement and the relevant definitions.
\end{proof}
\fi

Combined with \cref{prop:equivalence single communication class}, \cref{lem:closed classes} yields a convenient necessary and sufficient condition for convergence on~\(\classesmax\).
\begin{proposition}\label{prop:convergent on Xm iff regular on maxclasses}
An upper transition operator~$\utranop$ is convergent on $\classesmax$ if and only if \(\classmaxk{1}\), \dots, \(\classmaxk{M}\) are regular.
\end{proposition}
\begin{proof}
$\utranop$ is clearly convergent on $\classesmax$ if and only if it is convergent on all \(\classmaxk{1}\), \dots, \(\classmaxk{M}\).
The result now follows because, for all $k\in\inti{1}{M}$, 
$\classmaxk{k}$ is a closed communication class, so it follows from \cref{lem:closed classes} (for $\class=\classmaxk{k}$) and \cref{prop:equivalence single communication class} (applied to $\smash{\restrutranop[\classmaxk{k}]}$) that $\smash{\utranop}$ is convergent on $\classmaxk{k}$ if and only if $\classmaxk{k}$ is regular (for $\utranop$).
\end{proof}

\subsection{Maximal classes are absorbing}
With convergence on \(\classesmax\) completely covered, we now move to upper transition operators for which the closed class \(\classesmax\) is absorbing.
Recall from \cref{lem:linear then classesmax absorbing} that this is \emph{always} the case for linear transition operators.

More generally, we first look at convergence on an \emph{arbitrary} closed class~\(\class\) that is absorbing.
\begin{lemma}\label{lem:closed and absorbing convergence}
	Let $\class$ be an absorbing closed class for the upper transition operator~\(\utranop\).
	Then for any \(f\in\rfns\), if \(\gr{\utranop^n f}_{n\in\nats}\) converges on~$\class$, it converges (on $\stsp$) as well.
\end{lemma}
In our proof for this result, we'll also rely on \citep[Lemma~39]{2021TJoens}, which we repeat here for the sake of clarity.
\begin{lemma}
\label{lem:natan}
Consider an upper transition operator~\(\utranop\) with an absorbing closed class~\(\class\).
Then, for all \(\epsilon\in\posreals\), there is some \(n_\epsilon\in\nats\) such that for all \(n\geq n_\epsilon\) and all \(f\in\rfns\), \(\abs{\utranop^nf-\utranop^n(f\indica{\class})}\leq\max\abs{f}\epsilon\).
\end{lemma}
\begin{proof}[Proof of \cref{lem:closed and absorbing convergence}]
	It suffices to show that the limit set~\(\limitset=\set{\limit_1, \dots, \limit_{\per}}\) has period~\(\per=1\).

	For all \(k\in\inti{1}{\per}\), since \(\limit_k=\lim_{n\to+\infty}\utranop^{\per n}\limit_k\) by definition, it follows from \cref{lem:natan} that \(\limit_k=\lim_{n\to+\infty}\utranop^{\per n}(\limit_k\indica{\class})\). Since \(\restr{\limit_1}{\class}=\dots=\restr{\limit_{\per}}{\class}\) by assumption, and therefore also \(\limit_1\indica{\class}=\dots=\limit_{\per}\indica{\class}\), this implies that \(\limit_1=\cdots=\limit_{\per}\), proving that \(\per=1\).
\end{proof}
One immediate and interesting consequence of \cref{lem:closed and absorbing convergence} is that for an absorbing closed class~\(\class\), convergence on~\(\class\) is equivalent to convergence (on \(\stsp\)).
\begin{corollary}\label{cor:absorbing iff convergent}
Let $\class$ be a closed class that is absorbing.
Then $\utranop$ is convergent on~$\class$ if and only if it is convergent.
\end{corollary}
Of course, the prime example of such an absorbing closed class is the union of all maximal communication classes~\(\classesmax\) whenever \(\classtnl=\emptyset\).
Combined with \cref{prop:convergent on Xm iff regular on maxclasses}, this yields a necessary and sufficient condition for convergence for the case \(\classtnl=\emptyset\).

\begin{corollary}\label{cor:Xna empty}
	Consider any upper transition operator~\(\utranop\) such that \(\classtnl=\emptyset\).
	Then \(\utranop\) is convergent if and only if \(\classmaxk{1}\), \dots, \(\classmaxk{M}\) are regular.
\end{corollary}
For linear transition operators, for which we know from \cref{lem:linear then classesmax absorbing} that \(\classtnl=\emptyset\), this result specialises to the following necessary and sufficient condition for convergence.
\begin{corollary}
\label{cor:linear:iff}
    A linear transition operator~\(\tranop\) is convergent if and only if \(\classmaxk{1}\), \dots, \(\classmaxk{M}\) are regular.
\end{corollary}

\subsection{Maximal classes are not absorbing}
We have one more step to take: to allow for transient states from which \(\classesmax\) is not lower reachable, meaning that \(\classtnl\neq\emptyset\); as we've indicated before, this can only happen for non-linear transition operators.
As a first step, we establish that whenever \(\classtnl\neq\emptyset\), we can always restrict the upper transition operator~\(\smash{\utranop}\) to~\(\classtnl\).
\begin{lemma}
\label{lem:nonabs has nonempty restrcredset}
	For any upper transition operator~\(\utranop\) with \(\classtnl\neq\emptyset\), \(\restrutranop[\classtnl]\) is well defined.
\end{lemma}
\begin{proof}
	Let \(\gr{\credset{x}}_{x\in\stsp}\) be the family of credal sets corresponding to~\(\utranop\).
	We need to show that \(\restrcredset[\classtnl]{x}\neq\emptyset\) for all \(x\in\classtnl\).
	So fix any \(x\in\classtnl\).
	Let \(\gr{\class_n}_{n\in\nnegints}\) be as defined in \cref{lem:iterativeabs} with \(\class_0=\classesmax\), and recall from there that \(\classtnl=\stsp\setminus\class_n\) for all \(n\geq k\).
	Recall furthermore that
	\begin{equation*}
	\class_k=\class_{k+1}=\class_k\cup\{y\in\stsp\setminus\class_k\colon\ltranop\indica{\class_k}(y)>0\},
	\end{equation*}
	which, due to \ref{def:utranop:bounds}, implies that \(\ltranop\indica{\class_k}(x)=0\) because \(x\in\classtnl=\stsp\setminus\class_k\).
	Consider any \(\pmf\in\credset{x}\) such that \(\prev_{\pmf}(\indica{\class_k})=\min_{\altpmf\in\credset{x}} \prev_\altpmf(\indica{\class_k})=\ltranop \indica{\class_k} (x)=0\).
	Then clearly \(\pmf(y)=0\) for all \(y\in\class_k=\stsp\setminus\classtnl\), so \(\restr{\pmf}{\classtnl}\in\restrcredset[\classtnl]{x}\) by the definition in \cref{eqn:restrcredset:masses}, implying that \(\restrcredset[\classtnl]{x}\neq\emptyset\).
\end{proof}

Our next step is to use \(\smash{\restrutranop[\classtnl]}\) to decompose \(\classtnl\) into its own maximal and transient states, similarly to how we used \(\smash{\utranop}\) to decompose $\stsp$ into \(\classmaxk{1}\), \dots, \(\classmaxk{M}\), \(\classtl\) and \(\classtnl\). During this process, to avoid confusion, we will adopt \(\classmaxk{{\smash{1}}}^1\), \dots, \(\classmaxk{{\smash{M_1}}}^1\), \(\classtlsmash{1}\) and \(\classtnlsmash{1}\) as an alternative notation for \(\classmaxk{1}\), \dots, \(\classmaxk{M}\), \(\classtl\) and \(\classtnl\), respectively.
We now repeat the subdivision from before, but for \(\smash{\nonabsrestrutranop[1]}\) and \(\classtnlsmash{1}\) rather than \(\utranop\) and $\stsp$: we let \(\classmaxk{1}^2\), \dots, \(\classmaxk{M_2}^2\) denote the maximal classes of~\(\nonabsrestrutranop[1]\), let \(\classesmax^2\) denote the union of these maximal classes, let \(\smash{\classtlsmash{2}}\) collect all states in \(\classtnlsmash{1}\setminus\classesmax^2\) from which \(\classesmax^2\) is lower reachable by \(\nonabsrestrutranop[1]\), and let \(\classtnlsmash{2}\) collect all states in \(\classtnlsmash{1}\setminus\classesmax^2\) from which \(\classesmax^2\) isn't lower reachable by \(\smash{\nonabsrestrutranop[1]}\).
If \(\classtnlsmash{2}\neq\emptyset\), we repeat this process with \((\nonabsrestrutranop[1])_{\classtnl^2}\)---so the restriction of \(\nonabsrestrutranop[1]\) to~\(\smash{\classtnlsmash{2}}\)---to similarly decompose~\(\smash{\classtnlsmash{2}}\).
This notation is a bit unwieldy though. Luckily, the following lemma implies that we can use the simpler notation \(\nonabsrestrutranop[2]\) instead because \(\smash{\nonabsrestrutranop[2]=(\nonabsrestrutranop[1])_{\classtnl^2}}\).
\begin{lemma}
\label{lem:multiple cuts helper}
	Consider any upper transition operator~\(\utranop\) and any two classes \(\class,\altclass\subseteq\stsp\) such that \(\class\supseteq\altclass\).
	If \(\restrutranop[\class]\) and \(\gr[\big]{\restrutranop[\class]}_{\altclass}\) are well defined, then \(\gr[\big]{\restrutranop[\class]}_{\altclass}=\restrutranop[\altclass]\).
\end{lemma}
\ifarxiv
\begin{proof}
	Since \(\restrutranop[\class]\) and \(\gr[\big]{\restrutranop[\class]}_{\altclass}\) are well defined by assumption, let \((\restrcredset[\class]{x})_{x\in\class}\) and \(((\restrcredset[\class]{x})_{\altclass})_{x\in\class}\) be their corresponding families of credal sets.
	Fix some~\(x\in\class\).
	Because \(\stsp\supseteq\class\supseteq\altclass\) by the assumptions in the statement,
	\begin{align*}
	\big(\restrcredset[\class]{x}\big)_{\altclass}
	&= \set{\restr{\altpmf}{\altclass}\colon \altpmf\in\restrcredset[\class]{x}, \gr{\forall y\in\class\setminus\altclass}~\altpmf\gr{y}=0}\\
	&= \set{\restr{(\restr{\pmf}{\class})}{\altclass}\colon \pmf\in\credset{x}, \gr{\forall y\in\stsp\setminus\altclass}~\pmf\gr{y}=0}\\
	&= \set{\restr{\pmf}{\altclass}\colon \pmf\in\credset{x}, \gr{\forall y\in\stsp\setminus\altclass}~\pmf\gr{y}=0}\\
	&= \restrcredset[\altclass]{x},
	\end{align*}
	and therefore \(\gr[\big]{\restrutranop[\class]}_{\altclass}=\restrutranop[\altclass]\).
\end{proof}
\fi
We continue the process of repeated subdivisions until we reach a depth~\(d\in\nats\) such that \(\classtnlsmash{d}=\emptyset\)---this is always the case, as \(\stsp\) is finite.
Then by construction,
\(\classmaxk{{\smash{1}}}^1,\classmaxk{{\smash{2}}}^1, \dotsc, \classmaxk{{\smash{M_d}}}^d, \classtlsmash{1}, \dotsc, \classtlsmash{d}\)
\(\)
 partitions~\(\stsp\), as depicted in \cref{fig:cuts}.
Finally, we let \(\smash{\classesmax^*\coloneq\cup_{\ell=1}^d\cup_{k=1}^{M_\ell}\classmaxk{{\smash{k}}}^\ell}\) and \(\classtlsmash{*}\coloneq\cup_{\ell=1}^d \classtlsmash{\ell}\) and, for the sake of convenience, let \(\classtnlsmash{0}\coloneq\stsp\).

\begin{figure}
	\centering
	\begin{tikzpicture}[div/.style={minimum width=3.8em, minimum height=2em, inner sep=2pt}, partit/.style={draw, rounded corners, thick}, max/.style={div, partit, color=teal, fill=teal!5}, abson/.style={div, color=magenta}, absob/.style={partit, color=magenta, fill=magenta!5}, nonabson/.style={div, color=cyan}, nonabsob/.style={partit, color=cyan}]
		\node[max] (cmax11) {\(\classmaxk{1}^1\)};
		\node[max, right=.25cm of cmax11] (cmax12) {\(\classmaxk{2}^1\)};
		\node[max, right=.25cm of cmax12] (cmax13) {\(\classmaxk{3}^1\)};
		\node[right=.25cm of cmax13, color=teal] (cmax1d) {\(\cdots\)};
		\node[max, right=.25cm of cmax1d] (cmax1M) {\(\classmaxk{M_1}^1\)};
		\coordinate (c1) at ($(cmax12.south)!.5!(cmax13.south)$);
		\node[abson, below=.25cm of c1] (cabs1) {\(\classtl^1\)};
		\node[nonabson, below=.25cm of cmax1M] (cnabs1) {\(\classtnl^1\)\vphantom{\(\classmaxk{1}^2\)}};

		\coordinate[xshift=.25cm, yshift=-.5cm] (l2l) at (cmax11.south west|-cabs1.south east);
		\coordinate[xshift=-.25cm, yshift=-.5cm] (l2r) at (cmax1M.south east|-cabs1.south east);
		\node[max, anchor=north west] (cmax21) at (l2l) {\(\classmaxk{1}^2\)};
		\node[max, right=.25cm of cmax21] (cmax22) {\(\classmaxk{2}^2\)};
		\node[max, anchor=north east] (cmax2M) at (l2r) {\(\classmaxk{M_2}^2\)};
		\node[, color=teal] (cmax2d) at ($(cmax22)!.5!(cmax2M)$) {\(\cdots\)};
		\coordinate (c2) at ($(cmax21.south)!.5!(cmax2d.south)$);
		\node[abson, below=.25cm of c2] (cabs2) {\(\classtl^2\)};
		\node[nonabson, below=.25cm of cmax2M] (cnabs2) {\(\classtnl^2\)};

		\path (cabs2.south) -| node[below=.5cm] (dm1) {\(\vdots\)} ($(cmax21)!.5!(cmax2M)$);
		\node[div, color=teal, below = .5cm of dm1] (cmaxdd) {\(\cdots\)};
		\node[max, left=.25cm of cmaxdd] (cmaxd1) {\(\classmaxk{1}^d\)};
		\node[max, right=.25cm of cmaxdd] (cmaxdM) {\(\classmaxk{M_d}^d\)};
		\node[abson, below=.25cm of cmaxdd] (cabsd) {\(\classtl^d\)};

		\begin{scope}[on background layer]
			\draw[absob] (cabs1.north-|cmax11.west) rectangle (cabs1.south-|cmax1d.east);
			\draw[absob] (cabs2.north-|cmax21.west) rectangle (cabs2.south-|cmax2d.east);
			\draw[absob] (cabsd.north-|cmaxd1.west) rectangle (cabsd.south-|cmaxdM.east);
			\coordinate[yshift=-.25cm] (p1) at (cabsd.south);
			\coordinate[yshift=-.5cm] (p2) at (cabsd.south);
			\coordinate[yshift=-.75cm] (p3) at (cabsd.south);
			\coordinate[yshift=-1cm] (p4) at (cabsd.south);
			\draw[dotted, nonabsob] ($(dm1.south-|cmaxdM.east)+(.25, -.25)$) -- ++ (0, .5);
			\draw[dotted, nonabsob] ($(dm1.south-|cmaxd1.west)+(-.25, -.25)$) -- ++ (0, .5);
			\path[nonabsob] ($(dm1.south-|cmaxdM.east)+(.25, -.25)$) |- (p1) -| ($(dm1.south-|cmaxd1.west)+(-.25, -.25)$);
			\draw[nonabsob] (cnabs2.north east) |- (p2) -| ($(cabs2.south-|cmax21.west)+(0, -.25)$) -| (cnabs2.north west) -- cycle;
			\draw[nonabsob] (cnabs1.north east) |- (p3) -| ($(cabs1.south-|cmax11.west)+(0, -.25)$) -| (cnabs1.north west) -- cycle;
		\end{scope}
		\draw[partit] ($(cmax11.north west)+(-.25, .25)$) |- (p4) -| ($(cmax1M.north east)+(.25, .25)$) -- cycle;
	\end{tikzpicture}
	\caption{Modified Venn diagram of the state space~\(\stsp\)}
	\label{fig:cuts}
\end{figure}

This partitioning of the state space~\(\stsp\) allows us to present the following (general) sufficient condition for convergence of an upper transition operator.
\begin{theorem}
	\label{the: sufficient condition for convergence}
	Consider any upper transition operator~\(\utranop\).
	If \(\classmaxk{1}^\ell\), \dots, \(\classmaxk{M_\ell}^\ell\) are \(\nonabsrestrutranop\)-regular for all \(\ell\in\inti{1}{d}\), then \smash{\(\utranop\)} is convergent.
\end{theorem}
Before heading to its proof, let's apply this result to our running example.

\begin{runex}
\label{example:running example fifth part}
It's clear from \cref{fig:running example second part} that \(\smash{\restrutranop[\classtnl^1]}\) has a single maximal communication class that is given by \(\set{\std,\ste}\), so \(\classesmax^2=\classmaxk{1}^2=\set{\std,\ste}\) and \(\classtlsmash{2}=\classtnlsmash{2}=\emptyset\).
This finishes our decomposition of the state space.
To conclude that \(\utranop\) is convergent, let's check whether \(\classmaxk{1}^1\) and \(\classmaxk{2}^1\) are both \(\utranop\)-regular, and whether \(\classesmax^2\) is \(\smash{\nonabsrestrutranop[1]}\)-regular.
Since \(\classmaxk{1}^1\) and \(\classmaxk{2}^1\) both consist of a single element, it is immediate that they are \(\utranop\)-regular.
Furthermore, since the upper accessibility graph \(\uagraph[{\smash{\nonabsrestrutranop[1]}}]\) in \cref{fig:running example second part} has cyclicity \(1\), \(\classesmax^2\) is \(\nonabsrestrutranop[1]\)-regular.
\end{runex}

As a first step towards proving \cref{the: sufficient condition for convergence}, we present a result that enables us to use the \(\smash{\nonabsrestrutranop}\)-regularity of \(\classmaxk{k}^\ell\) to establish the convergence (and even ergodicity) of \(\utranop\) on \(\classmaxk{k}^\ell\).
\begin{lemma}
\label{lem:restriction converges then the limit is constant}
	Consider an upper transition operator~\(\utranop\) and a non-empty class~\(\class\) such that \(\smash{\restrutranop}\) is defined.
	If a maximal class of \(\restrutranop\) is \(\restrutranop\)-regular, then \(\utranop\) is ergodic on this class.
\end{lemma}
\begin{proof}
	Consider a maximal class \(\altclass\subseteq\class\) of \(\restrutranop\) that is \(\restrutranop\)-regular, fix any \(f\in\rfns\) and consider its limit set~\(\limitset=\{\limit_1, \dots, \limit_{\per}\}\).
	It suffices to show that there is a constant~\(\smash{\phi\in\rfns[\altclass]}\) such that \(\phi=\restr{\limit_1}{\altclass}=\cdots=\restr{\limit_{\per}}{\altclass}\).

	Fix any \(\smash{k\in\set{1, \dots, \per}}\). It's immediate from \cref{prop:limitdominates,lem:closed classes,lem:multiple cuts helper} that \(\gr{\gr{\restrutranop[\altclass]}^n \gr{\restr{\limit_k}{\altclass}}}_{n\in\nats}\) converges to some constant function, say \(\phi_k\in\rfns[\altclass]\), such that \(\phi_k\geq\min \restr{\limit_k}{\altclass}\) and, unless \(\restr{\limit_k}{\altclass}\) is a constant, \(\phi_k>\min\restr{\limit_k}{\altclass}\).
	Furthermore, due to \cref{eqn:inequality with restrutranop}, 
	\begin{equation*}
		\phi_k
		=\lim_{n\to+\infty}\gr{\restrutranop[\altclass]}^{\per n} \restr{\limit_k}{\altclass}
		\leq\lim_{n\to+\infty}\restr{\gr{\utranop^{\per n}\limit_k}}{\altclass}
		= \restr{\limit_k}{\altclass}.
	\end{equation*}
	It must therefore be that $\restr{\limit_k}{\altclass}$ is a constant, because otherwise $\min\restr{\limit_k}{\altclass}<\phi_k\leq\restr{\limit_k}{\altclass}$.
	Since $\min\restr{\limit_k}{\altclass}\leq\phi_k\leq\restr{\limit_k}{\altclass}$ and $\restr{\limit_k}{\altclass}$ and $\phi_k$ are constants, it follows that $\phi_k=\restr{\limit_k}{\altclass}$ is a constant function.

	It therefore suffices to show that \(\phi_1=\cdots=\phi_{\per}\).
	To this end, observe that for every \(\smash{k\in\set{1, \dots, \per}}\), it follows from \cref{eqn:inequality with restrutranop} that, with \(\phi_{\per+1}\coloneq\phi_1\),
	\begin{align*}
	\phi_k
	&=\lim_{n\to+\infty}\gr{\restrutranop[\altclass]}^{n+1} \restr{\limit_k}{\altclass} \\
	&=\lim_{n\to+\infty}\gr{\restrutranop[\altclass]}^n\restrutranop[\altclass] \restr{\limit_k}{\altclass} \\
	&\leq\limsup_{n\to+\infty}\gr{\restrutranop[\altclass]}^{n}\restr{\gr{\utranop\limit_k}}{\altclass} \\
	&=\lim_{n\to+\infty}\gr{\restrutranop[\altclass]}^n \restr{\limit_{k+1}}{\altclass}= \phi_{k+1},
	\end{align*}
	whence \(\phi_1\leq\phi_2\leq\dots\leq\phi_{\per}\leq\phi_1\), as required.
\end{proof}

If \(\classesmax^*\) were an absorbing closed class for $\utranop$, we could combine this result with \cref{lem:closed and absorbing convergence} to prove \cref{the: sufficient condition for convergence}; it is however not.
We therefore derive another upper transition operator~\(\utranop^*\) from~\(\utranop\) as follows: for all \(f\in\rfns\), let \(\utranop^*f(x)\coloneq\utranop f(x)\) for all \(x\in\classtlsmash{*}\) and let \(\utranop^*f(x)\coloneq\restrutranop[\classmaxk{k}^\ell](\restr{f}{\classmaxk{k}^\ell})(x)\) for all  \(\smash{x\in\classmaxk{k}^\ell}\), with \(\ell\in\inti{1}{d}\) and \(k\in\inti{1}{M_\ell}\);
it's immediate from \cref{lem:closed classes,lem:nonabs has nonempty restrcredset,lem:multiple cuts helper} that \(\utranop^*\) is well defined and it's easy to verify that it's an upper transition operator.
By \cref{eqn:inequality with restrutranop}, this definition ensures that
\begin{equation}
\label{eqn:utranop* vs utranop}
	\utranop^*f
	\leq \utranop f
	\quad\text{for all } f\in\rfns,
\end{equation}
with equality on~\(\classtlsmash{*}\).
That \(\utranop^*\) does satisfy the conditions in \cref{lem:closed and absorbing convergence} is our next result.
\begin{lemma}
\label{lem:classesmax is closed and absorbing}
Consider an upper transition operator~\(\utranop\).
Then \(\classesmax^*\) is an absorbing closed class for~\(\utranop^*\).
\end{lemma}
\begin{proof}
For all \(\ell\in\inti{1}{d}\) and \(k\in\inti{1}{M_\ell}\), \(\classmaxk{k}^\ell\) is a closed class for \(\utranop^*\) because, for all \(x\in\classmaxk{k}^\ell\) and \(y\in\stsp\setminus\classmaxk{k}^\ell\),
\(\utranop^*\indica{y}(x)=\restrutranop[\classmaxk{k}^\ell](0)(x)\overset{\smash[t]{\text{\ref{def:utranop:bounds}}}}{=}0\).
Since \(\classesmax^*\) is a union of such closed classes, it is itself closed as well.

	To show that \(\classesmax^*\) is absorbing for~\(\utranop^*\), we'll prove by induction that, for all $\ell\in\inti{1}{d}$, $\classesmax^*$ is lower reachable by~$\utranop^*$ from all states in $\classtlsmash{1:\ell}\coloneq\cup_{i=1}^\ell\classtlsmash{i}$.
	For $\ell=d$, we then find that $\classesmax^*$ is lower reachable by $\smash{\utranop^*}$ from all states in $\cup_{i=1}^d\classtlsmash{i}=\classtlsmash{*}=\stsp\setminus\classesmax^*$, or equivalently, that $\classesmax^*$ is indeed absorbing for~$\utranop^*$.

	For the base case $\ell=1$, we need to show that $\classesmax^*$ is lower reachable by $\smash{\utranop^*}$ from all states in $\classtlsmash{1}$.
	Consider any $x\in\classtlsmash{1}=\classtl$.
	Since $\classesmax$ is lower reachable by $\utranop$ from all states in $\classtl$, there is some $n\in\nats$ such that $\ltranop^n\indica{\classesmax}(x)>0$.
	Since $\classesmax\subseteq\classesmax^*$ and $\ltranop^*\geq\ltranop$ [\cref{eqn:utranop* vs utranop} and conjugacy], it follows that
	\begin{equation*}
		(\ltranop^*)^n\indica{\classesmax^*}(x)
		\geq
		\ltranop^n\indica{\classesmax^*}(x)
		\overset{\smash[t]{\text{\ref{def:utranop:monotonicity}}}}{\geq}
		\ltranop^n\indica{\classesmax}(x)>0,
	\end{equation*}
	whence $\classesmax^*$ is indeed lower reachable by~$\utranop^*$ from $x$.

	For the induction step, we assume that $\classesmax^*$ is lower reachable by $\utranop^*$ from all states in $\classtlsmash{1:\ell}$ for some $1\leq\ell<d$, and set out to prove that the same is then true for $\classtlsmash{1:\ell+1}=\classtlsmash{1:\ell}\cup\classtlsmash{\ell+1}$.

	By definition, $\classtlsmash{\ell+1}$ contains the states in $\classtnlsmash{\ell}\setminus\classesmax^{\ell+1}$ from which~$\classesmax^{\ell+1}$ is lower reachable by~$\nonabsrestrutranop[\ell]$. Hence, if we let $\smash{\class_0\coloneq\classesmax^{\ell+1}}$ and, for all \(n\in\nnegints\),
	\begin{equation*}
	\class_{n+1}
	\coloneq\class_n\cup\{x\in\classtnl^\ell\setminus\class_n\colon \nonabsrestrltranop[\ell]\indica{\class_n}(x)>0\},
	\end{equation*}
	then since $\classesmax^{\ell+1}$ is closed for~$\nonabsrestrutranop[\ell]$, \cref{lem:iterativeabs} says
	 that \(\smash{\classtlsmash{\ell+1}=\class_k\setminus\classesmax^{\ell+1}}\) for some \(k\in\nnegints\).

	We'll prove by induction over~\(n\) that $\class_n\setminus\classesmax^{\ell+1}$ is a set of states from which~$\classesmax^*$ is lower reachable by~$\utranop^*$.
	For $n=k$, we then find that $\classesmax^*$ is lower reachable by $\utranop^*$ from all states in $\classtlsmash{\ell+1}=\class_k\setminus\classesmax^{\ell+1}$.

	The base case $n=0$ is trivially true because $\class_0\setminus\classesmax^{\ell+1}=\emptyset$.
	For the induction step, we assume that, for some $n\in\inti{0}{k-1}$, $\classesmax^*$ is lower reachable from $\class_n\setminus\classesmax^{\ell+1}$ by~$\utranop^*$, and set out to prove that $\classesmax^*$ is then lower reachable from $\class_{n+1}\setminus\class_n$ by~$\utranop^*$ as well.

	To that end, consider any $x\in\class_{n+1}\setminus\class_n$.
	Since $\classesmax^*$ is a closed class for~$\utranop^*$, we know from \cref{lem:iterativeabs} that there is some~$K\in\nats$ such that $(\ltranop^*)^K\indica{\classesmax^*}(y)>0$ for all $y\in\classesmax^*$ and all $y\in\stsp\setminus\classesmax^*$ from which~$\classesmax^*$ is lower reachable by~$\utranop^*$.
	We now set out to prove that $(\ltranop^*)^{K+1}\indica{\classesmax^*}(x)>0$, thereby indeed establishing that $\classesmax^*$ is lower reachable from $x$ by $\utranop^*$, as required.

	Since it follows from our induction hypothesises that $\classtlsmash{1:\ell}$ and $\class_n\setminus\classesmax^{\ell+1}$ are sets from which~$\classesmax^*$ is lower reachable by~$\utranop^*$, we know that $(\ltranop^*)^K\indica{\classesmax^*}(y)>0$ for all $y\in\gr{\stsp\setminus\classtnl^\ell}\cup C_n$.
	Since $(\ltranop^*)^K\indica{\classesmax^*}\geq0$ by~\ref{def:utranop:bounds}, this implies that there is some $\alpha>0$ such that $(\ltranop^*)^K\indica{\classesmax^*}\geq\alpha\indica{\gr{\stsp\setminus\classtnl^\ell}\cup \class_n}$.
	Hence,
	\begin{align}
	(\ltranop^*)^{K+1}\indica{\classesmax^*}(x)
	&=\ltranop^*\big((\ltranop^*)^{K}\indica{\classesmax^*}\big)(x)\notag\\
	\overset{\smash[t]{\text{\ref{def:utranop:monotonicity}}}}&{\geq}\ltranop^*\big(\alpha\indica{\gr{\stsp\setminus\classtnl^\ell}\cup C_n}\big)(x)\notag\\
	\overset{\smash[t]{\text{\ref{def:utranop:positively homogeneous}}}}&{=}\alpha\ltranop^*\big(\indica{\gr{\stsp\setminus\classtnl^\ell}\cup C_n}\big)(x)\notag\\
	&=\alpha\ltranop\big(\indica{\gr{\stsp\setminus\classtnl^\ell}\cup C_n}\big)(x),\label{eq:inductionabsorbed}
	\end{align}
	where the last equality follows from the definition of~\(\utranop^*\) because $x\in\class_{n+1}\setminus\class_n\subseteq\class_k\setminus\classesmax^{\ell+1}=\classtlsmash{\ell+1}\subseteq\classtlsmash{*}$.

	Let \((\credset{y})_{y\in\stsp}\) be the family of credal sets that correspond with \(\utranop\).
	Since \(\credset{x}\) is closed by definition, there is some $p\in\credset{x}$ such that
	\begin{align*}
	\ltranop\big(\indica{\gr{\stsp\setminus\classtnl^\ell}\cup C_n}\big)(x)
	&=\prev_p(\indica{\gr{\stsp\setminus\classtnl^\ell}\cup C_n})\\
	&=\prev_p(\indica{\stsp\setminus\classtnl^\ell})+\prev_p(\indica{C_n}),
	\end{align*}
	where for the second equality we used that \(\smash{\gr{\stsp\setminus\classtnl^\ell}\cap \class_n=\emptyset}\), which holds because \(\class_n\subseteq\classtnl^\ell\) by definition.
	Since both terms in the sum are clearly non-negative, it suffices to establish that at least one of them is positive to show that $\ltranop^*\big(\indica{\gr{\stsp\setminus\classtnl^\ell}\cup C_n}\big)(x)>0$, which then implies $\smash{(\ltranop^*)^{K+1}\indica{\classesmax^*}(x)>0}$ due to \cref{eq:inductionabsorbed}.
	To that end, we'll assume that $p(y)=0$ for all $y\in\stsp\setminus\classtnl^\ell$, or equivalently, that $\smash{\prev_p(\indica{\stsp\setminus\classtnl^\ell})=0}$ and show that this implies that $\prev_p(\indica{C_n})>0$.

	Since $p(y)=0$ for all $y\in\stsp\setminus\classtnl^\ell$, \cref{eqn:restrcredset:masses} tells us that $\restr{\pmf}{\classtnl^\ell}\in\restrpmfset[\classtnl^{\smash{\ell}}]{x}$.
	This implies that
	\begin{equation*}
	\prev_p(\indica{C_n})=\prev_{\restr{\pmf}{\classtnl^\ell}}(\indica{C_n})\geq\nonabsrestrltranop[\ell]\indica{\class_n}(x)>0,
	\end{equation*}
	where for the final inequality we use the fact that $x\in C_{n+1}\setminus C_n$.
\end{proof}

Finally, then, we can lay out our proof for \cref{the: sufficient condition for convergence}.

\begin{proof}[Proof of \cref{the: sufficient condition for convergence}]
	Fix \(f\in\rfns\) and consider its limit set~\(\limitset=\set{\limit_1, \dots, \limit_{\per}}\). It suffices to show that~\(\per=1\).

	For any \(\ell\in\inti{1}{d}\) and \(k\in\inti{1}{M_\ell}\), since \(\classmaxk{k}^\ell\) is a maximal class of \(\nonabsrestrutranop\) that is \(\nonabsrestrutranop\)-regular, it's immediate from \cref{lem:restriction converges then the limit is constant} that \(\utranop\) is ergodic on \(\smash{\classmaxk{k}^\ell}\), which implies that \((\utranop^n\limit_1)_{n\in\nats}\) converges to a constant on \(\classmaxk{k}^\ell\), and hence, \(\restr{\limit_1}{\classmaxk{k}^\ell}=\cdots=\restr{\limit_{\per}}{\classmaxk{k}^\ell}\) is constant.
	Consequently, for all \(i\in\inti{1}{\per}\),
	\begin{align*}
	\limit_{i+1}(x)=\limit_i(x)
	\overset{\smash[t]{\text{\ref{def:utranop:bounds}}}}{=}
	\restrutranop[\classmaxk{k}^\ell]\gr{\restr{\limit_i}{\classmaxk{k}^\ell}}(x)
	=\utranop^*\limit_i(x)
	\end{align*}
	for all \(x\in\classmaxk{k}^\ell\). Since also
\(\limit_{i+1}(x)=\utranop \limit_i(x)=\utranop^*\limit_i(x)\) for all \(i\in\inti{1}{\per}\) and $x\in\classtnlsmash{*}$,
we conclude that \(\smash{\limit_{i+1}
=\utranop^*\limit_i}\) for all \(i\in\inti{1}{\per}\).
Consequently, the limit set of~\(\gr{\gr{\utranop^*}^n \limit_1}_{n\in\nats}\) is \(\set{\limit_1, \dots, \limit_{\per}}\).

Since \(\restr{\limit_1}{\classmaxk{k}^\ell}=\cdots=\restr{\limit_{\per}}{\classmaxk{k}^\ell}\) for all \(\ell\in\inti{1}{d}\) and \(k\in\inti{1}{M_\ell}\), this implies that \(\smash{\gr{\gr{\utranop^*}^n \limit_1}_{n\in\nats}}\)  converges on~\(\classesmax^*\).
	Now recall from \cref{lem:classesmax is closed and absorbing} that for \(\utranop^*\), \(\classesmax^*\) is an absorbing closed class.
	Consequently, it follows from \cref{lem:closed and absorbing convergence} that \(\gr{\gr{\utranop^*}^n\limit_1}_{n\in\nats}\) converges, so \(\limit_1=\dots=\limit_{\per}\) and therefore \(\per=1\).
\end{proof}

\section{The finitely generated case}
A natural follow-up question is whether the sufficient condition in \cref{the: sufficient condition for convergence} is also a necessary one.
This is the case, at least if the upper transition operator~\(\smash{\utranop}\) is \emph{finitely generated}, meaning that
it is compatible with a
family~\(\gr{\pmfset{x}}_{x\in\stsp}\) of \emph{finite} sets of pmfs.\footnote{In fact, the reader may want to verify that it suffices for~\(\utranop\) to be \emph{sufficiently finitely generated}, meaning that it is compatible with a
family~\(\gr{\pmfset{x}}_{x\in\stsp}\) of sets of pmfs such that for all \(x\in\classtnl\), \(\pmfset{x}\) is a finite set.}
\begin{proposition}
\label{prop:convergence for finitely generated:nec alt}
	Consider a finitely generated upper transition operator~\(\utranop\).
	If \(\utranop\) is convergent, then for all \(\ell\in\set{1, \dots, d}\), \(\classmaxk{1}^\ell\), \dots, \(\classmaxk{M_\ell}^\ell\) are \(\nonabsrestrutranop\)-regular.
\end{proposition}
Our proof follows relatively straightforward from some intermediary results. Since it is rather instructive, we run through it in the main text.
First, recall from \cref{prop:convergent on Xm iff regular on maxclasses} that as \(\utranop\) is convergent, the maximal communication classes~\(\classmaxk{k}^1\) are regular for \(\nonabsrestrutranop[\smash{0}]=\utranop\).
So if \(d=1\), we're already done.
If on the other hand \(d>1\)---and therefore \(\classtnl\neq\emptyset\)---we turn our attention to the behaviour on~\(\classtnl\).
\begin{proposition}
\label{prop:convergence for finitely generated:nec part II}
	Consider a finitely generated upper transition operator~\(\utranop\).
	If \(\utranop\) is convergent and \(\classtnl\neq\emptyset\), then \(\restrutranop[\classtnl]\) is convergent.
\end{proposition}
\ifarxiv
Our proof for \cref{prop:convergence for finitely generated:nec part II} hinges on the following lemma.
\begin{lemma}
\label{lem:hannes helper}
	Consider a finitely generated upper transition operator~\(\utranop\) with \(\classtnl\neq\emptyset\), and fix some compatible family~\(\gr{\pmfset{x}}_{x\in\stsp}\) of finite sets of pmfs.
	Let
	\begin{align*}
		\kappa
		&\coloneq \max\set{\utranop\indica{\classtnl}(x)\colon x\in\classtl}\cup\{0\}
		< 1
	\intertext{and, with \(\mathcal{Q}\coloneq\bigcup_{x\in\classtnl}\set{\pmf\in\pmfset{x}\colon \prev_\pmf(\indica{\classtnl})<1}\),}
		\lambda
		&\coloneq \max \set{\prev_\pmf(\indica{\classtnl})\colon \pmf\in\mathcal{Q}}
		<1.
	\end{align*}
	Take any \(h\in\rfns[\classtnl]\) such that \(0\leq h\leq 1\), constants~\(\alpha,\beta\in\nposreals\) such that \(\alpha\leq\beta\) and a function~\(f\in\rfns\) such that \(\restr{f}{\classesmax}=\alpha\), \(\restr{f}{\classtl}\leq\beta\) and \(\restr{f}{\classtnl}=h\).
	Then \(\restr{(\utranop f)}{\classesmax}=\alpha\), \(\restr{(\utranop f)}{\classtl}\leq (1-\kappa)\beta+\kappa\) and, if \(\beta<1/(\lambda-1)\), \(\restr{(\utranop f)}{\classtnl}=\restrutranop[\classtnl] h\).
\end{lemma}
In its turn, the proof of \cref{lem:hannes helper} relies on the following two lemmas.
\begin{lemma}
\label{lem:utranop indica classtnl}
	For any upper transition operator~\(\utranop\),
	\begin{equation*}
		\utranop\indica{\classtnl}(x)
		< 1
		\quad\text{for all } x\in\classesmax\cup\classtl.
	\end{equation*}
\end{lemma}
\begin{proof}
	Since \(\classesmax\) is a closed class, we see that for all \(x\in\classesmax\),
	\begin{equation*}
		\utranop\indica{\classtnl}(x)
		\overset{\smash[t]{\text{\ref{def:utranop:subadditive}}}}{\leq} \sum_{y\in\classtnl} \utranop\indica{y}(x)
		= 0
		< 1.
	\end{equation*}

	From \cref{lem:iterativeabs} with \(\class=\classesmax\), we know that
	\begin{equation*}
		\classtnl
		= \stsp\setminus\class_k
		\subsetneq \stsp\setminus\class_{k-1}
		\subsetneq \cdots
		\subsetneq \stsp\setminus\class_0
		= \stsp\setminus\classesmax,
	\end{equation*}
	with the increasing sequence \(\class_0, \dots, \class_k\) as defined there.
	Hence, for all \(n\in\inti{0}{k}\), \(\indica{\classtnl}\leq\indica{\stsp\setminus \class_n}\) and therefore
	\begin{equation*}
		\utranop \indica{\classtnl}
		\overset{\smash[t]{\text{\ref{def:utranop:monotonicity}}}}{\leq} \utranop\indica{\stsp\setminus \class_n}.
	\end{equation*}
	For all \(x\in\classtl\), it follows from the definition in \cref{lem:iterativeabs} that there is some \(n\in\inti{0}{k-1}\) such that \(\ltranop\indica{\class_n}(x)>0\), or equivalently [use \ref{prop:utranop:indicators}], \(\utranop\indica{\stsp\setminus\class_n}(x)<1\). It follows from this and the preceding inequality that \(\utranop \indica{\classtnl}(x)<1\).
\end{proof}
\begin{lemma}
\label{lem:lambda is well defined}
	Consider a finitely generated upper transition operator~\(\utranop\) with \(\classtnl\neq\emptyset\), and fix some compatible family~\(\gr{\pmfset{x}}_{x\in\stsp}\) of finite sets of pmfs.
	Then
	\begin{equation*}
		\bigcup_{x\in\classtnl} \{\pmf\in\pmfset{x}\colon \prev_\pmf(\indica{\classtnl})<1\}
		\neq \emptyset.
	\end{equation*}
\end{lemma}
\begin{proof}
Assume \emph{ex absurdo} that \(\prev_\pmf(\indica{\classtnl})=1\) for all \(\pmf\in\pmfset{x}\) and \(x\in\classtnl\).
	Then for all \(x\in\classtnl\),
	\begin{align*}
		\utranop\indica{\stsp\setminus\classtnl}(x)
		&= \max\{\prev_\pmf(\indica{\stsp\setminus\classtnl})\colon \pmf\in\pmfset{x}\} \\
		&= \max\{\prev_\pmf(1-\indica{\classtnl})\colon \pmf\in\pmfset{x}\} \\
		&= 1-\min\{\prev_\pmf(\indica{\classtnl})\colon \pmf\in\pmfset{x}\} = 0.
	\end{align*}
	But this means that \(\classtnl\) is a closed class, so that it must include at least one maximal communication class, which of course contradicts the fact that \(\classtnl\subseteq\stsp\setminus\classesmax\) only includes non-maximal states.
\end{proof}
\begin{proof}[Proof of \cref{lem:hannes helper}]
	That \(\kappa\) and \(\lambda\) are well defined and satisfy the given inequalities follows from \cref{lem:utranop indica classtnl,lem:lambda is well defined} and the fact that \(\pmfset{x}\) is finite for all \(x\in\classtnl\).

	For any maximal communication class~\(\classmaxk{k}\), it follows from \cref{lem:closed classes} that
	\begin{equation*}
		\restr{(\utranop f)}{\classmaxk{k}}
		= \restrutranop[\classmaxk{k}] (\restr{f}{\classmaxk{k}})
		= \restrutranop[\classmaxk{k}] \alpha
		\overset{\smash[t]{\text{\ref{def:utranop:bounds}}}}{=} \alpha.
	\end{equation*}
	Consequently, \(\restr{(\utranop f)}{\classesmax}=\alpha\).
	Furthermore, by the conditions on~\(f\) in the statement,
	\begin{multline}
		f
		\leq \beta \indica{\stsp\setminus\classtnl} + (\max h) \indica{\classtnl} \\
		\leq \beta \indica{\stsp\setminus\classtnl} + \indica{\classtnl}
		= \beta + (1-\beta) \indica{\classtnl},
		\label{eqn:proof of helper lemma:f and indica classtnl}
	\end{multline}
	and therefore
	\begin{equation*}
		\utranop f
		\overset{\smash[t]{\text{\ref{def:utranop:monotonicity}}}}{\leq} \utranop (\beta + (1-\beta) \indica{\classtnl})
		\overset{\smash[t]{\text{\ref{def:utranop:constant additivity},\ref{def:utranop:positively homogeneous}}}}{=} \beta + (1-\beta) \utranop\indica{\classtnl}.
	\end{equation*}
	It follows from this and the definition of \(\kappa\) that
	\begin{equation*}
		\restr{(\utranop f)}{\classtl}
		\leq \beta + (1-\beta) \kappa
		= (1-\kappa)\beta + \kappa.
	\end{equation*}

	To complete our proof, we assume that \(\beta<1/(\lambda-1)\) [which implies that \(\beta+(1-\beta)\lambda<0\) since \(\lambda<1\)] and set out to show that then \(\restr{(\utranop f)}{\classtnl}=\restrutranop[\classtnl] h\), where $\restrutranop[\classtnl]$ is well-defined because of \cref{lem:nonabs has nonempty restrcredset}.

	To this end, fix some \(x\in\classtnl\). We need to show that \(\utranop f(x)=\restrutranop[\classtnl] h(x)\).
	Note that for all \(\pmf\in\pmfset{x}\),
	\begin{equation*}
		\gr[\big]{\gr{\forall y\in\stsp\setminus\classtnl}~\pmf\gr{y}=0}
		\Leftrightarrow \prev_{\pmf}\gr{\indica{\classtnl}}=1.
	\end{equation*}
	From this, the definition of \(\restrutranop[\classtnl]\) and the fact that \(\restr{f}{\classtnl}=h\), we infer that
	\begin{align*}
		\restrutranop[\classtnl] h (x)
		&= \max\set[\big]{\prev_{\altpmf}(h)\colon \altpmf\in\restrpmfset[\classtnl]{x}} \\
		&= \max\set[\big]{\prev_{\restr{\pmf}{\classtnl}}(h)\colon \pmf\in\pmfset{x}, \prev_{\pmf}\gr{\indica{\classtnl}}=1} \\
		&= \max\set[\big]{\prev_{\pmf}(f)\colon \pmf\in\pmfset{x}, \prev_{\pmf}\gr{\indica{\classtnl}}=1}.
	\end{align*}
	On the other hand,
	\begin{align*}
		\utranop f(x)
		&= \max\{\prev_\pmf(f)\colon \pmf\in\pmfset{x}\} \\
		&= \max\{\prev_\pmf(f)\colon \pmf\in\pmfset{x}, \prev_\pmf(\indica{\classtnl})=1\}\\&\qquad\qquad\cup\{\prev_\pmf(f)\colon \pmf\in\pmfset{x}, \prev_\pmf(\indica{\classtnl})<1\}.
	\end{align*}
	Consequently, it suffices for us to prove that
	\begin{equation*}
		\max\{\prev_\pmf(f)\colon \pmf\in\pmfset{x}, \prev_\pmf(\indica{\classtnl})<1\}
		\leq \restrutranop[\classtnl] h (x).
	\end{equation*}
	Now for any \(\pmf\in\pmfset{x}\), it follows from \cref{eqn:proof of helper lemma:f and indica classtnl} and the monotonicity and linearity of~\(\prev_\pmf\) that
	\begin{equation*}
		\prev_\pmf(f)
		\leq \beta + (1-\beta) \prev_\pmf(\indica{\classtnl}).
	\end{equation*}
	From this and the conditions in the statement, it follows that
	\begin{multline*}
		\max\{\prev_\pmf(f)\colon \pmf\in\pmfset{x}, \prev_\pmf(\indica{\classtnl})<1\} \\
		\leq \beta + (1-\beta)\lambda
		< 0.
	\end{multline*}
	Since we know from \ref{def:utranop:bounds} that \(\restrutranop[\classtnl] h(x)\geq\min h\geq 0\), this finalises our proof.
\end{proof}
Finally, then, everything is in place for us to prove \cref{prop:convergence for finitely generated:nec part II}.
\begin{proof}[Proof of \cref{prop:convergence for finitely generated:nec part II}]
	Let \(\class\coloneq\classtnl\).
	Suppose towards contradiction that there is some \(g\in\rfns[\class]\) such that \((\gr{\restrutranop}^n g)_{n\in\nats}\) does not converge.
	We'll now use this \(g\) to construct some \(f\in\rfns\) for which
\(\restr{(\utranop^n f)}{\class}=\gr{\restrutranop}^n g\) for all \(n\in\nats\), which is a contradiction because \(\utranop\) was assumed to be convergent.

	Because \(\restrutranop\) satisfies \ref{def:utranop:bounds}, \(g\) must be non-constant, so \(\max g>\min g\).
	Furthermore, since \(\restrutranop\) satisfies \ref{def:utranop:constant additivity} and \ref{def:utranop:positively homogeneous}, we may assume without loss of generality that \(0\leq g\leq 1\)---otherwise, continue with \(g'\coloneq(g-\min g)/(\max g-\min g)\).

	From \cref{lem:iterativeabs} for \(\classesmax\), combined with \ref{def:utranop:bounds}, it follows that there are some \(\epsilon>0\) and \(k\in\nats\) such that
	\begin{equation}
	\label{eqn:proof of fgn:epsilon}
		\ltranop^n\indica{\classesmax}(x)
		\geq\epsilon
		\quad\text{for all } n\geq k, x\in\classtl.
	\end{equation}

	Since \(\utranop\) is assumed to be finitely generated, it's induced by some family~\(\gr{\pmfset{x}}_{x\in\stsp}\) such that \(\pmfset{x}\) is finite for all \(x\in\stsp\).
	Let \(\kappa\) and \(\lambda\) be as defined in \cref{lem:hannes helper}.
	Then because \(0\leq\kappa<1\), \(0\leq\lambda<1\) and \(\epsilon>0\), we can fix some \(b_0\in\reals_{<0}\) such that
	\begin{align}
	\label{eqn:condition on beta0:one}
	b_0&<\left(\frac1{\lambda-1}-1\right)\frac1\epsilon+1
	\shortintertext{and}
	\label{eqn:condition on beta0:two}
	b_0&<\left(\frac1{\lambda-1}-1\right)\frac1{(1-\kappa)^{k-1}}+1.
	\end{align}
	For all \(n\in\nats\), let \(b_n\coloneq(1-\kappa)^n(b_0-1)+1\); it's easy to verify that in this way,
	\begin{equation}
	\label{eqn:bn recursive}
		b_n
		= (1-\kappa)b_{n-1}+\kappa
		\quad\text{for all } n\in\nats;
	\end{equation}
	that is, $b_n$ is a convex mixture of $b_{n-1}$ and $1$, with mixture coefficient \(0\leq\kappa<1\).
	Since $b_0<0$, this implies that $b_0\leq b_1\leq\cdots\leq b_{k-1}$.
	It therefore follows from the definition of~\(b_{k-1}\) and \cref{eqn:condition on beta0:two} that
	\begin{equation}
	\label{eqn:bn is increasing}
		b_0
		\leq b_1
		\leq \cdots
		\leq b_{k-1}
		< \frac{1}{\lambda-1}.
	\end{equation}

	Finally, we let
	\begin{equation*}
		f
		\colon \stsp\to\reals\colon x\mapsto
		\begin{cases}
			g(x) &\text{if }x\in\classtnl \\
			b_0 &\text{otherwise} \\
		\end{cases}
	\end{equation*}
	and set out to verify that \(\restr{(\utranop^n f)}{\class}=\gr{\restrutranop}^n g\) for all \(n\in\nats\).

	First, we show by induction that for all \(n\in\set{1, \dots, k}\), \(\restr{(\utranop^n f)}{\classesmax}=b_0\), \(\restr{(\utranop^n f)}{\classtl}\leq b_n\) and \(0\leq\restr{(\utranop^n f)}{\class}=(\restrutranop)^n g\leq 1\).
	For the base case \(n=1\), this follows from \cref{lem:hannes helper}---with \(\alpha=b_0\), \(\beta=b_0\) and \(h=g\)---due to our definition of~\(f\) and \cref{eqn:bn recursive,eqn:bn is increasing}.
	The inductive step is similar: if these three properties hold for \(n\in\set{1, \dots, k-1}\), then it follows from \cref{lem:hannes helper}---with \(\utranop^n f\) here in the role of \(f\) there, \(\alpha=b_0\), \(\beta=b_n\) and \(h=(\restrutranop)^n g\)---, \cref{eqn:bn recursive,eqn:bn is increasing} and the boundedness of \((\restrutranop)^n g\) [\ref{def:utranop:bounds}] that they also hold for \(n+1\).

	Second, we show by induction that for all \(n\geq k\), \(\restr{(\utranop^n f)}{\classesmax}=b_0\) and \(0\leq\restr{(\utranop^n f)}{\class}=(\restrutranop)^n g\leq 1\).
	The base case \(n=k\) was already shown in the previous paragraph.
	Hence, for the inductive step we assume that these two properties hold for some \(n\geq k\).
	As \(f\leq b_0\indica{\classesmax}+(1-\indica{\classesmax})\) by construction, it follows that for all \(x\in\classtl\),
	\begin{align*}
		\utranop^n f(x)
		\overset{\smash[t]{\text{\ref{def:utranop:monotonicity}}}}&{\leq} \utranop^n [b_0\indica{\classesmax}+(1-\indica{\classesmax})](x) \\
		&= (b_0-1)\ltranop^n \indica{\classesmax}(x)+1 \\
		\overset{\smash[t]{\text{\eqref{eqn:proof of fgn:epsilon}}}}&{\leq} (b_0-1)\epsilon+1.
	\end{align*}
	Because of this, the induction hypothesis and \cref{eqn:condition on beta0:one}, it follows from \cref{lem:hannes helper}---with \(\utranop^n f\) here in the role of \(f\) there, \(\alpha=b_0\), \(\beta=(b_0-1)\epsilon+1\) and \(h=\gr{\restrutranop}^n g\)---that \(\restr{(\utranop^{n+1} f)}{\class}=b_0\) and \(0\leq\restr{(\utranop^{n+1} f)}{\class}=(\restrutranop)^{n+1} g\leq 1\), as required.
\end{proof}
\fi
Since \(\nonabsrestrutranop[1]=\restrutranop[\classtnl]\) is convergent, we may again use \cref{prop:convergent on Xm iff regular on maxclasses} to infer that \(\classmaxk{1}^2\), \dots, \(\classmaxk{M_2}^2\) are regular for \(\nonabsrestrutranop[1]\).
Now if \(d=2\), we're done.
If on the other hand \(d>2\), we want to repeat the argument, and for this we need that \(\nonabsrestrutranop[1]\) is finitely generated. This is however clearly the case. Indeed, as \(\smash{\utranop}\) is finitely generated, it is compatible with a family~\(\gr{\pmfset{x}}_{x\in\stsp}\) of sets of pmfs that are all finite (and therefore closed).
Since restrictions of finite credal sets are finite themselves, \cref{lem:restrutranop via closed family} therefore implies that \(\nonabsrestrutranop[1]\) is finitely generated as well.

Since \(\nonabsrestrutranop[{\smash{1}}]\) is also convergent, \cref{prop:convergence for finitely generated:nec part II} tells us that \(\smash{\nonabsrestrutranop[{\smash{2}}]
}\) is convergent, and then \cref{prop:convergent on Xm iff regular on maxclasses} establishes that \(\classmaxk{1}^3\), \dots, \(\classmaxk{M_3}^3\) are \(\nonabsrestrutranop[2]\)-regular.

Repeated application of the same argument until depth~\(d\)---with \(\classtnl^d=\emptyset\)---eventually results in (a proof for) \cref{prop:convergence for finitely generated:nec alt}.

To conclude our treatment for finitely generated upper transition operators, we combine \cref{the: sufficient condition for convergence,prop:convergence for finitely generated:nec alt,prop:convergence for finitely generated:nec part II} into the following strong result.
\begin{theorem}
\label{the:necessary and sufficient for finitely generated}
	For a finitely generated upper transition operator~\(\utranop\), the following are equivalent:
	\begin{enumerate}[label=\upshape(\roman*)]
		\item \(\utranop\) is convergent;
		\item the maximal communication classes~\(\classmaxk{1}\), \dots, \(\classmaxk{M}\) are regular and \(\restrutranop[\classtnl]\) is convergent;
		\item for all \(\ell\in\set{1, \dots, d}\), the maximal communication classes~\(\classmaxk{1}^\ell\), \dots, \(\classmaxk{M_\ell}^\ell\) are \(\nonabsrestrutranop\)-regular.
	\end{enumerate}
\end{theorem}

\Cref{the:necessary and sufficient for finitely generated} immediately leads to an algorithm to determine whether a finitely generated upper transition operator is convergent.
It consists in recursively taking the following steps, starting from \(\ell=1\): (i) construct the upper accessibility graph~\(\uagraph[{\nonabsrestrutranop[\ell-1]}]\), (ii) determine the maximal communication classes~\(\smash{\classmaxk{1}^\ell}\), \dots, \(\smash{\classmaxk{M_\ell}^\ell}\) and their cyclicity using one of the standard algorithms; (iii) determine \(\classtnl^\ell\) with the recursive procedure in \cref{lem:iterativeabs}; and (iv) if \(\smash{\classtnl^\ell\neq\emptyset}\) determine~\(\nonabsrestrutranop[\ell]\) [via finite sets of pmfs thanks to~\cref{lem:restrutranop via closed family,lem:multiple cuts helper}], increment~\(\ell\) and repeat, otherwise stop.

\section{Wrapping things up}

The results above give rise to another follow-up question: can we generalise \cref{prop:convergence for finitely generated:nec part II} (and then also \cref{prop:convergence for finitely generated:nec alt,the:necessary and sufficient for finitely generated}) to the general case of upper transition operators that need not be finitely generated?
As is clear from \cref{prop:convergent on Xm iff regular on maxclasses}, it's always necessary that the maximal communication classes are regular.
Unfortunately, though, for an upper transition operator that is \emph{not} finitely generated, it's no longer necessary for convergence that \(\restrutranop[\classtnl]\) is convergent, making it impossible to generalise \cref{prop:convergence for finitely generated:nec part II} to this case.
What follows is a straightforward example of a convergent upper transition operator~\(\smash{\utranop}\) that is not finitely generated such that \(\restrutranop[\classtnl]\) is not convergent---because its sole communication class  is not regular.
\subsection{Counterexample}
\label{ssec:counterexample}
Let \(\stsp\coloneq\{\sta,\stb,\stc\}\) and \(\pmf_\epsilon\coloneqq \epsilon^2\indica{\sta}+\epsilon\indica{\stb}+(1-\epsilon-\epsilon^2)\indica{\stc}\) for all \(\epsilon\in\cci{0}{1}\).
Then the sets of transition pmfs
\begin{equation*}
\pmfset{\sta}\coloneq\set{\indica{\sta}},
\pmfset{\stb}\coloneq\set{\indica{\sta}}\cup\set{\pmf_\epsilon\colon \epsilon\in\cci{0}{\nicefrac12}},
\pmfset{\stc}\coloneq\set{\indica{\sta},\indica{\stb}}
\end{equation*}
induce the upper transition operator~\(\utranop\) given for all \(f\in\rfns\) and \(x\in\stsp\) by
\begin{equation*}
\utranop f(x)
= \begin{cases}
f(\sta) &\text{if }x=\sta, \\
\max\set{f(\sta)}\cup\set{\prev_{\pmf_\epsilon}\gr{f}\colon \epsilon\in\cci{0}{\nicefrac12}} &\text{if }x=\stb, \\
\max\set*{f(\sta),f(\stb)} &\text{if }x=\stc.
\end{cases}
\end{equation*}
From the upper accessibility graph~\(\uagraph\) depicted in \cref{fig:counterexample}, it's clear that the upper transition operator~\(\utranop\) has one maximal class: \(\classesmax^1=\classmaxk{1}^1=\set{\sta}\).
\begin{figure}
\centering
	\subfloat[{\(\uagraph[\utranop]\)}]{%
		\label{fig:counterexample}
		\begin{tikzpicture}[node distance=1.5cm,on grid,auto]
			\node[state] (a) {\(\sta\)};
			\node[state] (b) [below left=of a] {\(\stb\)};
			\node[state] (c) [below right=of a] {\(\stc\)};
			\path[->,>={Stealth[round]}]
			(b) edge[bend left]  (a)
			(b) edge[loop left]  (b)
			(b) edge[bend right] (c)
			(a) edge[loop above] (a)
			(c) edge[bend right] (b)
			(c) edge[bend right] (a)
			(c) edge [loop right, opacity=0] (c);
		\end{tikzpicture}
	}%
	\subfloat[{\(\uagraph[{\nonabsrestrutranop[1]}]\)}]{
		\label{fig:counterexample cut}%
		\begin{tikzpicture}[node distance=1.5cm, on grid, auto]
			\node[state] (b) {\(\stb\)};
			\node[state] (c) [right=of b] {\(\stc\)};
			\path[->,>={Stealth[round]}]
			(b) edge[bend right] (c)
			(c) edge[bend right] (b);
		\end{tikzpicture}
	}
\caption{}
\end{figure}
Since \ref{prop:utranop:indicators} implies that $\smash{\ltranop\indica{\sta}=1-\utranop\indica{\set{\stb,\stc}}=1-\indica{\set{\stb,\stc}}=\indica{\sta}}$,
it follows that
\begin{align*}
\classtl^1
=
\set{x\in\set{\stb,\stc}\colon (\exists n\in\nats)~\ltranop^n\indica{\sta}(x)>0}
=
\emptyset,
\end{align*}
so \(\classtnl^1=\set{\stb,\stc}\).

Since \(\pmfset{\stb}\) and \(\pmfset{\stc}\) are closed, \cref{lem:restrutranop via closed family} implies that \(\nonabsrestrutranop[1]\) is compatible with the restricted sets of pmfs
\begin{equation*}
\smash{\restrpmfset[\classtnl^{\smash{1}}]{\stb}=\set{\indica{\stc}}}
\text{~~and~~}
\smash{\restrpmfset[\classtnl^{\smash{1}}]{\stc}=\set{\indica{\stb}}};
\end{equation*} it is therefore easy to see that for all \(g\in\rfns[\classtnl^1]\),
\begin{equation*}
	\nonabsrestrutranop[1]g\gr{\stb}
	= g(\stc)
	\quad\text{and}\quad
	\nonabsrestrutranop[1]g\gr{\stc}
	= g(\stb).
\end{equation*}
Its upper accessibility graph~\(\uagraph[{\nonabsrestrutranop[{\smash{1}}]}]\) is depicted in \cref{fig:counterexample cut}. It has only one (and therefore maximal) communication class \(\classesmax^2=\set{\stb,\stc}=\classtnl^1\), so \(\classtl^2=\emptyset=\classtnl^2\).
Since it's clear that \(\classesmax^2\) has cyclicity \(2\), it is not regular and therefore, due to \cref{prop:equivalence single communication class}, \(\nonabsrestrutranop[1]\) is not convergent.
Consequently, the sufficient condition for convergence in \cref{the: sufficient condition for convergence} is \emph{not} met, nor is the necessary condition in \cref{prop:convergence for finitely generated:nec part II}.

Nonetheless, as we will now show, \(\utranop\) is convergent. In particular, for all \(f\in\rfns\) and \(x\in\stsp\),
\begin{equation*}
	\lim_{n\to+\infty} \utranop^n f(x)
	= \begin{cases}
		f(\sta) &\text{if }x=\sta \\
		\max f &\text{if } x\neq\sta.
	\end{cases}
\end{equation*}
Since $\utranop f\gr{\sta}=f\gr{\sta}=\min\utranop f$ and $\max\utranop f=\max f$, we can assume without loss of generality that $\min f=f\gr{\sta}$.
It follows immediately from the expression for~\(\utranop\) that \(\utranop^n f\gr{\sta}=f\gr{\sta}\) for all \(n\in\nats\), so we can focus on the value in~\(\stb\) and \(\stc\).

If \(f(\sta)\leq f(\stb)=f(\stc)=\max f\), then it follows from the expression for~\(\utranop\) that \(\utranop^n f=f\) for all \(n\in\nats\).
Hence, it remains for us to look at the limiting behaviour on~\(\stb\) and \(\stc\) in case \(f\gr{\sta}<\max f\) with $f(\stb)\neq f(\stc)$.

For all \(n\in\nats\), let \(f_n^-\coloneq\min\{\utranop^n f\gr{\stb},\utranop^n f\gr{\stc}\}\) and \(f_n^+\coloneq\max\{\utranop^n f\gr{\stb},\utranop^n f\gr{\stc}\}\).
We continue with a general observation: for any \(h\in\rfns\), we write \(h=\gr{h\gr{\sta}, h\gr{\stb}, h\gr{\stc}}\) and observe that
\begin{align}
h(\sta)\leq h(\stb)<h(\stc)
&\Rightarrow
\utranop h
= \gr{h(\sta), h(\stc), h(\stb)};
\label{eq:switch}\\
h(\sta)\leq h(\stc)<h(\stb)
&\Rightarrow
\begin{cases}\utranop h(\sta)=h(\sta),\\  h(\stc)<\utranop h(\stb)\leq h(\stb), \\ \utranop h(\stc)= h(\stb).\end{cases}
\label{eq:strict}
\end{align}
From \cref{eq:switch,eq:strict} and the definition of \(\utranop\), it follows that (i) \(f_n^+=\max f\) for all \(n\in\nats\); and (ii) the sequence \(\gr{f_n^-}_{n\in\nats}\) is non-decreasing in \(\cci{f(\sta)}{\max f}\), and therefore converges to a limit $\lambda^-$ with \(f(\sta)\leq\lambda^-\leq\max f\).
We now need to show that \(\lambda^-=\max f\), so we assume towards contradiction that \(f(\sta)\leq\lambda^-<\max f\).
On the one hand, the orbit \(\gr{\utranop^{n}f}_{n\in\nats}\) then has a limit set \(\limitset=\set{\limit_1, \limit_2}\) of period $\per=2$, which alternates between \(\limit_1=h_{\stc}\coloneq\gr{f\gr{\sta}, \max f, \lambda^-}\) and \(\limit_2=h_{\stb}\coloneq\gr{f\gr{\sta}, \lambda^-, \max f}\). This implies in particular that $h_{\stb}=\limit_2=\utranop\limit_1=\utranop h_{\stc}$.
On the other hand, it follows from \cref{eq:strict} that \(\utranop h_{\stc}(\stb) > h_{\stc}(\stc)=h_{\stb}(\stb)\); the contradiction we were after!

\subsection{Outlook}
In future work, we'd like to find out whether for the general, not necessarily finitely generated case, there is some---necessarily other---equivalent characterisation of convergence that can also be easily checked.
Our preliminary research has already revealed that it is indeed possible to come up with an interesting equivalent condition, but we've not yet succeeded at translating this condition into one that can be easily verified in practice.

We also plan to scrutinise the relation between our \cref{the: sufficient condition for convergence} and \citeauthor{2003Akian}'s~\citep{2003Akian} Theorem~5.5, which says that \(\utranop\) is convergent if all of the strongly connected components of their `critical graph'~\(\cgraph\) have cyclicity~\(1\), where \(\cgraph\) is defined in terms of the accessibility graphs of all the linear transition operators $\tranop$ that are dominated by $\utranop$.
While we strongly believe our condition is more convenient to verify than theirs, there seems to be a strong connection between our condition and theirs, and we wouldn't be surprised if they turn out to be equivalent.

\additionalinfo

\begin{acknowledgements}
	A lot of the ideas and some of the results---for example the counterexample in \cref{ssec:counterexample}, \cref{prop:equivalence single communication class,prop:convergent on Xm iff regular on maxclasses,cor:Xna empty,prop:convergence for finitely generated:nec part II}---in this contribution have their roots in Hanne Asselman's Master's thesis \citep{2024Asselman}, which was supervised by the authors.
	These ideas were further developed by the authors and Natan T'Joens during a first research retreat and by the authors during a second research retreat.
	The authors thank Hanne and Natan for their valuable input, as well as for their blessing to continue this line of research without them.

	Our thanks also goes out to the three anonymous reviewers for their constructive comments.

	Jasper De Bock's research is supported by Research Foundation – Flanders (FWO), project number G0ACO25N.
	Alexander Erreygers is a post-doctoral research fellow of Ghent University's Special Research Fund; his research was also supported by Research Foundation – Flanders (FWO), project number 3G028919.
	Floris Persiau is currently a post-doctoral fellow of the Belgian American Educational Foundation; his research was also supported by Research Foundation – Flanders (FWO), project number 11H5523N.

\end{acknowledgements}


\begin{authorcontributions}
The order in which the authors are listed is alphabetical; it is not intended to reflect the extent of their contribution. All authors have contributed equally.
\end{authorcontributions}

\printbibliography

\end{document}